\documentclass{amsart}
\usepackage{amsxtra}
\usepackage{amssymb,amsmath,amsthm}
\usepackage{color}
\usepackage[all]{xy}
\usepackage{hyperref}
\usepackage{appendix}
\theoremstyle{plain}

\newtheorem{theorem}{Theorem}[section]
\newtheorem{lemma}[theorem]{Lemma}
\newtheorem{proposition}[theorem]{Proposition}
\newtheorem{corollary}[theorem]{Corollary}
\newtheorem{definition}[theorem]{Definition}
\theoremstyle{definition}
\newtheorem{example}[theorem]{Example}
\newtheorem{remark}[theorem]{Remark}

\newcommand \lb[1]{\label{#1}}

\newcommand{\rmap}       {\longrightarrow}

\newcommand{\id}         {{\mathrm {Id}}}

\newcommand{\pr}         {{\mathrm{pr}}}

\newcommand{\Id}         {\mathrm{Id}}

\newcommand{\lie}[1]{\mathfrak{#1}}

\newcommand{\source}   {\mathsf{s}}
\newcommand{\target}   {{\mathsf{t}}}

\newcommand{\tg}       {\mathsf t}
\newcommand{\s}        {\mathsf s}
\newcommand{\core}        {\mathsf c}
\newcommand{\rr}       {\rightrightarrows}

\newcommand{\V}{{\mathcal{V}}}

\newcommand{\R}{\mathbb{R}}

\newcommand{\inv}{^{-1}}

\newcommand{\mx}{\mathfrak{X}}
\newcommand{\dr}{\mathbf{d}}

\newcommand{\ldr}[1]{{{\pounds}}_{#1}}
\newcommand{\ip}[1]{{\mathbf{i}}_{#1}}
\newcommand{\an}[1]{\arrowvert_{#1}}

\DeclareMathOperator{\dom}{Dom}
\DeclareMathOperator{\erz}{span}
\DeclareMathOperator{\Exp}{Exp}

\begin{document}
\title{Foliated groupoids and their infinitesimal data}


\author{M. Jotz}
\address{Georg-August-Universit\"at G\"ottingen\\
Mathematisches Institut\\
Bunsenstr. 3-5\\
37073 G\"ottingen\\
Germany}
\thanks{The first author was partially supported by Swiss NSF grant 200021-121512 and by the
\emph{Dorothea Schl\"ozer Programme} of the  University of G\"ottingen.}
\email{mjotz@uni-math.gwdg.de}


\author{C. Ortiz}
\address{Departamento de Matem\'atica\\
Universidade Federal do Paran\'a\\
Setor de Ciencias Exatas - Centro Politecnico\\
81531-990 Curitiba\\
Brazil}
\email{cristian.ortiz@ufpr.br}


\begin{abstract}

In this work, we study  Lie groupoids equipped with multiplicative foliations and the corresponding 
infinitesimal data. 
We determine the infinitesimal counterpart of a multiplicative foliation in terms of its core and sides 
together with a partial connection satisfying special properties, giving rise to the concept of 
 IM-foliation on a Lie algebroid. The main result of this paper shows that if $G$ is a 
source simply connected Lie groupoid with Lie algebroid $A$, then there exists a one-to-one 
correspondence between multiplicative foliations on $G$ and IM-foliations on the Lie algebroid $A$.

\end{abstract}
\maketitle

\tableofcontents

\section{Introduction}\lb{intro}
This paper is the first part of a series of two articles devoted to the study 
of the infinitesimal data of Dirac groupoids, i.e. 
Lie groupoids equipped with Dirac structures suitably compatible with the groupoid multiplication. 

Since Dirac structures generalize Poisson bivectors, 
closed $2$-forms and regular foliations, multiplicative Dirac structures \cite{Ortiz08t,thesis}
provide a unified framework for the study of a variety of geometric structures compatible with a 
group(oid) multiplication, including Poisson groupoids and multiplicative closed $2$-forms. 
It is well known that Poisson groupoids are infinitesimally described by Lie bialgebroids \cite{MaXu00}.
 Similarly, the infinitesimal counterpart of multiplicative closed $2$-forms on a Lie groupoid are 
IM-$2$-forms on the corresponding Lie algebroid \cite{BuCrWeZh04,BuCaOr09}.
 Here we are concerned with Lie groupoids 
equipped with a multiplicative foliation. In particular, our main goal is 
to describe multiplicative foliations at the infinitesimal level.

Consider a Lie group $G$ with Lie algebra $\lie g$ and multiplication map 
$\mathsf m:G\times G\to G$. Then the tangent space $TG$ of $G$ is also a Lie group with 
multiplication map $T\mathsf m:TG\times TG\to TG$. It was observed in \cite{Jotz11b, Jotz11a, Ortiz08} 
that a multiplicative 
distribution $F_G\subseteq TG$, that is,  a distribution on $G$ that is  a subgroup
of $TG$, is the bi-invariant image of an ideal $\lie f \subseteq \lie g$.
Hence, it is  automatically an \emph{involutive subbundle} of $TG$ which is \emph{completely 
determined by its fiber $\lie f$ over the unit} of the Lie group.

If $G\rr M$ is a Lie groupoid, then the application of the tangent functor to all the groupoid maps 
gives rise to a Lie groupoid $TG\rr TM$, called the tangent groupoid. A distribution  
$F_G\subseteq TG$ on $G$ is said to be \textbf{multiplicative} if $F_G$ is a subgroupoid of 
$TG\rr TM$ over $F_M\subseteq TM$. In this more general situation it is no longer true that $F_G$ is automatically
involutive and has constant rank. Indeed, every smooth manifold can be viewed as a Lie groupoid 
over itself, and, in this case, any distribution $F\subseteq TM$ is multiplicative.

\medskip

We study in this paper  how the easy properties of multiplicative distributions on Lie groups 
can be extended to the more general framework of Lie groupoids. A \textbf{multiplicative foliation} 
on a Lie groupoid $G\rr M$ is a multiplicative \emph{subbundle} $F_G\subseteq TG$ which is also 
involutive with respect to the Lie bracket of vector fields. We say
that the pair $(G\rr M,F_G)$ is a 
\textbf{foliated groupoid}. The main goal of this paper is to describe foliated groupoids 
infinitesimally, that is, in terms of Lie algebroid data. 

First, we observe that given a Lie groupoid $G\rr M$ with Lie algebroid $A$, then the space 
of units $F_M\subseteq TM$ of a multiplicative foliation $F_G\subseteq TG$ is necessarily an involutive subbundle, and similarly the core $F_\core:=F_G\cap T^\s_MG$ of $F_G$ is a Lie subalgebroid 
of the Lie algebroid $A= T^\s_MG$. In particular, a multiplicative foliation $F_G\subseteq TG$ can be thought of as a 
$\mathcal{LA}$-subgroupoid of $TG\rr TM$. As a consequence, following \cite{Mackenzie00, Ortiz08t} 
we prove that the Lie algebroid of $F_G$ is, modulo the canonical identification
$A(TG)\simeq TA$, a subalgebroid $F_A\rmap F_M$ of the tangent 
Lie algebroid $TA\rmap TM$ which at the same time is an involutive subbundle $F_A\subseteq TA$. 
In this case we say that $(A,F_A)$ is a \textbf{foliated algebroid} and $F_A$ is referred to as 
a \textbf{morphic foliation} on $A$. We show that if $G\rr M$ is a source simply connected Lie 
groupoid with Lie algebroid $A$, then there is a one-to-one correspondence between multiplicative 
foliations on $G$ and morphic foliations on $A$. 

We prove that if $G$ is a Lie groupoid with 
Lie algebroid $A$ endowed with a multiplicative foliation $F_G$ with core $F_\core\subseteq A$ and units 
$F_M\subseteq TM$, then there is a well-defined partial
$F_M$-connection $\nabla$ on $A/F_\core$, which 
is naturally induced by the Bott connection associated to the foliation $F_G$ on $G$. We show that this 
connection is flat and allows us to compute explicitly the Lie algebroid of $F_G\rr F_M$ in 
terms of restrictions to $F_M\subseteq TM$ of special sections of the Lie algebroid 
$A(TG)\rmap TM$ of the tangent groupoid.
Furthermore, we show that the involutivity of the multiplicative subbundle $F_G\subseteq TG$ is completely encoded in some data involving only 
$F_M$, $F_\core$ and 
the flatness of the connection $\nabla$.

Given a foliated algebroid $(A,F_A)$, we show that the core 
$F_\core$ and side $F_M$ of $F_A$ are Lie subalgebroids of $A$ and $TM$, 
respectively. We prove that, here also, there is a flat partial $F_M$-connection on 
$A/F_\core$ which is naturally induced by
 the Bott connection defined by the foliation $F_A$ on $A$.  
In addition, we show that if $A$ integrates to a Lie groupoid $G$ and $F_A\subseteq TA$ 
integrates to a multiplicative foliation $F_G\subseteq TG$, then the partial connections 
induced by $F_A$ and $F_G$ coincide.
We say that 
the data $(A,F_M,F_\core,\nabla)$ is an \textbf{IM-foliation} on $A$.

\medskip

The main result of this work states that the data $(A,F_M,F_\core,\nabla)$ 
determines completely the foliated algebroid $(A,F_A)$ in the sense that $F_A$ 
can be reconstructed out of $(A,F_M,F_\core,\nabla)$.  In particular, 
if $G\rr M$ is a source simply connected Lie groupoid with Lie algebroid $A$, then 
there is a one-to-one correspondence between multiplicative foliations on $G$ and 
IM-foliations $(A,F_M,F_\core,\nabla)$   on $A$.

The partial connection induced by a multiplicative foliation already appeared in  \cite{Jotz11b},
where the existence of its parallel sections is used  to show that under some 
completeness conditions, there is a groupoid structure on the leaf space 
$G/F_G\rr M/F_M$ of the multiplicative foliation.
Analogously, we study the leaf space of a 
foliated algebroid $(A,F_A)$. We prove that, whenever the leaf space
$M/F_M$ is a 
smooth manifold and the partial connection has trivial holonomy,
 then there is a natural Lie algebroid structure on $A/F_A\rmap M/F_M$. We 
also show that if a foliated groupoid $(G\rr M,F_G)$ is such that the completeness conditions of 
\cite{Jotz11b} are fulfilled, then the Lie groupoid $G/F_G\rr M/F_M$ has Lie algebroid 
$A/F_A\rmap M/F_M$. Several examples of foliated groupoids and algebroids are 
discussed, including actions by groupoid and algebroid automorphisms, complex Lie groupoids and complex 
Lie algebroids 
and integrable regular Dirac structures versus regular presymplectic groupoids\footnote{Here, \emph{regular}
means that the characteristic distributions have constant rank.}.

\medskip

Given a Lie groupoid $G\rr M$ with Lie algebroid $A$, it was proved in \cite{CoDaWe87} 
that the cotangent bundle $T^*G$ has a Lie groupoid structure over $A^*$. In particular, the direct 
sum vector bundle $TG\oplus T^*G$ inherits a Lie groupoid structure over $TM\oplus A^*$. It 
is easy to see that if $F_G\subseteq TG$ is a multiplicative foliation on $G$, then its 
annihilator $F_G^\circ\subseteq T^*G$ is a Lie subgroupoid. Hence, the subbundle 
$F_G\oplus F_G^\circ\subseteq TG\oplus T^*G$ is an example of a multiplicative 
Dirac structure on $G$ (see \cite{Jotz11c, Ortiz08t, thesis}). Every multiplicative Dirac structure
induces Lie algebroid structures on the core and the base \cite{Jotz11c}. Furthermore
there is a Courant algebroid that is canonically associated to the Dirac groupoid \cite{Jotz11c}, 
which generalizes the Lie bialgebroid of a 
Poisson groupoid and is closely related to the partial connections studied in this paper.

The understanding of  multiplicative Dirac structures at the infinitesimal level
in the spirit of this paper
is a very interesting goal, because this
 generalizes simultaneously the correspondences between 
Poisson groupoids and Lie bialgebroids \cite{MaXu00}, 
multiplicative closed $2$-forms and IM-$2$-forms \cite{BuCrWeZh04,BuCaOr09},
and the results in this paper.
This infinitesimal description of Dirac groupoids 
can be made with similar techniques as here, by using the core and 
side algebroids of a multiplicative Dirac structure $\mathsf D_G$ on $G$ together with the Courant algebroid 
 associated to $\mathsf D_G$.
We will treat this problem in a separate paper \cite{DrJoOr11}.

\vspace*{0.5cm}

\paragraph{\textbf{Remark}}
Our result on the correspondence between IM-foliations 
and foliated algebroids (foliated groupoids) was already 
stated in \cite{Hawkins08}, but we were not aware of that while working on this subject.

\vspace*{0.5cm}

\paragraph{\textbf{Structure of the paper and main results}}
 In Section \ref{backgrounds}, we start by recalling background on Lie groupoid
theory, especially the tangent and cotangent groupoids associated to a Lie groupoid, as well
as the Lie algebroid structures on $A(TG)\to TM$ and on 
$TA\to TM$ in terms of linear and core sections.
We also discuss briefly the definition and some properties of flat partial connections.

In Section \ref{fol_gpds}, we define foliated Lie groupoids, together with the
associated subalgebroids $F_M\subseteq TM$ and $F_\core\subseteq A$.  Our first main results are summarized 
in the following theorem.
\vspace*{0.5cm}

\noindent{\bf Theorem 1} {\it
Let $G\rr M$ be a Lie groupoid endowed with a multiplicative subbundle $F_G\subseteq TG$.
Then there is a map
\begin{align*}
\tilde\nabla:\Gamma(F_M)\times\Gamma\left(A\right)&\to \Gamma\left(A/F_\core\right)
\end{align*}
such that the following holds:
$F_G$ is involutive if and only if
\begin{enumerate}
\item[a)] $F_M\subseteq TM$ is involutive,
\item[b)] for any $\bar X\in\Gamma(F_M)$, $\tilde\nabla(\bar X,\cdot)$ vanishes on sections of $F_\core$ and
\item[c)]  the induced map
\begin{align*}
\nabla:\Gamma(F_M)\times\Gamma\left(A/F_\core\right)&\to \Gamma\left(A/F_\core\right).
\end{align*}
is then a flat partial $F_M$-connection on $ A/F_\core$. 
\end{enumerate}
The connection $\nabla$ has in this case the following properties:
\begin{enumerate}
\item If $a\in\Gamma(A)$ is $\nabla$-parallel, i.e., 
$\nabla_{\bar X}\bar a=0$ for all $\bar X\in\Gamma(F_M)$, where $\bar
a$ is the class of $a$ in $\Gamma(A/F_\core)$, then $[a,b]\in\Gamma(F_\core)$ 
for all $b\in\Gamma(F_\core)$.
\item If $a,b\in\Gamma(A)$ are $\nabla$-parallel, then $[a,b]$ is also $\nabla$-parallel.
\item If $a\in\Gamma(A)$ is  $\nabla$-parallel, then $[\rho(a),\bar X]\in\Gamma(F_M)$ for all
$\bar X\in\Gamma(F_M)$. That is, $\rho(a)$ is $\nabla^{F_M}$-parallel, where $\nabla^{F_M}$ denotes the Bott connection defined by $F_M$.
\end{enumerate}
Furthermore,  $a\in\Gamma(A)$ is  $\nabla$-parallel if and only if the right invariant vector field $a^r$ on $G$ is parallel with respect to the Bott connection $\nabla^{F_G}$ defined by $F_G$. That is \[\left[a^r,\Gamma(F_G)\right]\subseteq \Gamma(F_G).\]}

In Section \ref{fol_alg}, we introduce and study the concept of foliated Lie algebroids. We compute explicitly  
the Lie algebroid of a multiplicative foliation in terms of parallel sections of the flat partial 
connection of Theorem 1, and we show the following integration result. 

\vspace*{0.5cm}
\noindent{\bf Theorem 2} {\it
Let $(G\rr M, F_G)$ be a foliated groupoid. Then 
$(A, F_A=j_G\inv(A(F_G)))$ is a foliated algebroid.\footnote{ Here,
$j_G:TA\to A(TG)$ is the canonical identification.}

Conversely, let $(A,F_A)$ be a foliated Lie algebroid. Assume that
$A$ integrates to a source simply connected Lie groupoid $G\rr M$. Then there is a unique multiplicative 
foliation $F_G$ on $G$ such that $F_A=j_G\inv(A(F_G))$.
}
\vspace*{0.5cm}

A foliated algebroid $(A,F_A)$ induces similar data $(A,F_M,F_\core,\nabla)$  as does
a foliated groupoid. 

\vspace*{0.5cm}
\noindent{\bf Theorem 3} {\it
Let $(A, F_A)$ be a foliated algebroid. Then $F_M\subseteq TM$
and $F_\core\subseteq A$ are subalgebroids and there 
is a natural partial $F_M$-connection $\nabla^A$ on $A/F^c$  induced by $F_A$.
This connection has the same properties as the connection found in Theorem 1.

The parallel sections of $\nabla^A$ are exactly
the sections $a\in\Gamma(A)$ such that  $[a^\uparrow,\Gamma_A(F_A)]\subseteq \Gamma_A(F_A)$, that is the core vector field $a^\uparrow$ on $A$ is parallel with respect to the Bott connection defined by $F_A$.}
\vspace*{0.5cm}

The data induced by a foliated groupoid and its foliated algebroid are the same, according to the next result.

\vspace*{0.5cm}
\noindent{\bf Theorem 4} {\it
Assume that $A$ integrates to a Lie groupoid $G\rr M$ and that $j_G(F_A)$ integrates to $F_G\subseteq TG$.
Then $\nabla=\nabla^{A}$.
}
\vspace*{.5cm}

The next theorem, which is proved in Section \ref{connection-as-infinitesimal-data},
 is the main result of this work.  We prove that the data $(A,F_M,F_\core,\nabla)$,
which is called an \textbf{IM-foliation} on $A$,
 is what should be considered as the infinitesimal data of a multiplicative foliation.

\vspace*{0.5cm}
\noindent{\bf Theorem 5} {\it
Let $(A\to M,\rho,[\cdot\,,\cdot])$ be a Lie algebroid
and $(A,F_M,F_\core,\nabla)$ an IM-foliation on $A$.
Then there exists a unique  morphic
 foliation $F_A$ on $A$  with side $F_M$ and core $F_\core$
such that the induced connection $\nabla^{A}$ is equal to $\nabla$.

If $A$ integrates to a source simply connected Lie groupoid $G\rr M$, there 
is a unique multiplicative foliation $F_G$ on $G$ 
with associated infinitesimal data $(A,F_M,F_\core,\nabla)$.
}
\vspace*{.5cm}

In Section \ref{quotient_algebroid} we study the leaf space of a foliated algebroid, and relate
it in the integrable case to the leaf space of the
corresponding foliated groupoid \cite{Jotz11b}.

\vspace*{0.5cm}
\noindent{\bf Theorem 6} {\it
Let $(A,F_M,F_\core,\nabla)$ be an IM-foliation on $A$ and
assume that the leaf space $M/F_M$ is a  smooth manifold and $\nabla$ has trivial holonomy.
Then there is natural Lie algebroid structure on $A/F_A$ over $M/F_M$ such that the projection 
$(\pi,\pi_M)$ is a Lie algebroid morphism.
\begin{align*}
\begin{xy}
\xymatrix{
A\ar[d]_{q_A}\ar[r]^{\pi}&A/F_A\ar[d]^{[q_A]}\\
M\ar[r]_{\pi_M}&M/F_M
}
\end{xy}
\end{align*}

If $(A,F_A)$ integrates to $(G\rr M,F_G)$ and $G/F_G\rr M/F_M$ is a Lie groupoid,
then $A/F_A$ with the structure above is the Lie algebroid of  $G/F_G\rr M/F_M$.
}
\vspace*{.5cm}

We show in the last section that regular Dirac manifolds 
provide an interesting example of IM-foliations. The corresponding 
multiplicative foliation is the kernel of the presymplectic groupoid integrating the Dirac structure.
Finally, we discuss shortly the relation of this work with the foliated algebroids 
in the sense of Vaisman \cite{Vaisman10b}.
\vspace*{0.5cm}

\paragraph{\textbf{Notation}}
Let $M$ be a smooth manifold. We will denote by $\mx(M)$ and $\Omega^1(M)$ the
spaces of (local) smooth sections of the tangent and the cotangent bundle,
respectively. For an arbitrary vector bundle $E\to M$, the space of
(local) sections of $E$ will be written $\Gamma(E)$. 
We will write $\dom(\sigma)$ for the open subset of the smooth manifold $M$
where the local section $\sigma\in\Gamma(E)$ is defined.

The flow of a vector field $X$ will be written $\phi^X_\cdot$,
 unless specified otherwise.

Let $f:M\to N$ be a smooth map between two smooth manifolds $M$ and $N$.
Then two vector fields $X\in\mx(M)$ and $Y\in\mx(N)$ are said to be \textbf{$f$-related} 
if $Tf\circ X=Y\circ f$ on $\dom(X)\cap f\inv(\dom(Y))$.
We write then $X\sim_f Y$.

The pullback or restriction of a vector bundle $E\to M$ to an embedded
 submanifold $N$ of $M$  will be written $E\an{N}$. In the special
case of the tangent and cotangent spaces of $M$, we will write $T_NM$ and $T^*_NM$.
If $f:M\to N$ is a smooth surjective submersion, we write
$T^fM$ for the kernel of $Tf:TM\to TN$.

\vspace*{0.5cm}

\paragraph{\textbf{Acknowledgements}}
We  would like to thank Henrique Bursztyn,  Thiago 
Drummond, Mathieu Sti\'enon  and Ping Xu 
for interesting discussions. 
We are especially thankful to Maria Amelia Salazar who showed us
a counterexample to a former version of Theorem 3.10. Thanks to her comments,
we improved considerably Section 3.3.

The first author wishes to 
thank Penn State University, IMPA and Universidade Federal do Paran\'a, Curitiba, where
parts of this work have been done. The second author thanks IMPA, 
EPFL and ESI for the hospitality during several stages of this project.

\vspace*{0.5cm}

\section{Background}\label{backgrounds}
\subsection{Basics on Lie Theory}

\subsubsection{Lie groupoids}

A \textbf{groupoid} $G$ over the \textbf{units} $M$ will be written $G\rr M$.
The \textbf{source} and \textbf{target} maps are denoted by
 $\source,\target:G\rmap M$ respectively,  the \textbf{unit section} $\epsilon:M\rmap G$, the
 inversion map 
$\mathsf i:G\rmap G$ and the  multiplication $\mathsf m:G_{(2)}\rmap G$, where 
$G_{(2)}=\{(g,h)\in G\times G\mid \target(h)=\source(g)\}$ is the set of composable 
groupoid pairs.  

A groupoid $G$ over $M$ is called a \textbf{Lie groupoid} if both $G$ and $M$ are smooth Hausdorff manifolds, 
the source and target maps $\source,\target:G\rmap M$ are surjective submersions, and all the other 
structural maps are smooth. 
Throughout this work we only consider Lie groupoids.

Let $G$ and $G'$ be Lie groupoids over $M$ and $M'$, respectively. A \textbf{morphism} of Lie 
groupoids is a smooth map $\Phi:G'\rmap G$ over $\phi:M'\to M$
that is compatible with all the structural maps. When $\Phi$ is injective we say that 
$G'$ is a \textbf{Lie subgroupoid} of $G$. See \cite{Mackenzie05} for more details. 

\subsubsection{Lie algebroids}
 A \textbf{Lie algebroid} is a vector bundle $q_A:A\rmap M$ equipped with a Lie bracket 
$[\cdot,\cdot]_A$ on the space of smooth sections $\Gamma_M(A)$ and a vector bundle map 
$\rho_A:A\rmap TM$ called the \textbf{anchor}, such that  

$$[a,fb]_A=f[a,b]_A+(\ldr{\rho_A(a)}f)b,$$

\noindent for every $a,b\in\Gamma_M(A)$ and $f\in C^{\infty}(M)$.

Let $A$ and $A'$ be Lie algebroids over $M$ and $M'$, respectively. A \textbf{Lie algebroid morphism} 
is a bundle map $\phi:A'\rmap A$ over $\phi:M'\to M$
which is compatible with the anchor and bracket, see \cite{Mackenzie05} for more details. 
If $\phi:A'\rmap A$ is injective, we say that $A'$ is a \textbf{Lie subalgebroid} of $A$.

\subsubsection{The Lie functor}
Let $G\rr M$ be a Lie groupoid. The Lie algebroid of $G$ is defined in this paper
to be $AG=T^\s_MG$, with anchor $\rho_{AG}:=T\tg\an{AG}$
and bracket $[\cdot\,,\cdot]_{AG}$ defined by using right invariant vector fields.

We will write $A(\cdot)$
for the functor that sends Lie groupoids to Lie algebroids 
and Lie groupoid morphisms  to Lie algebroid
morphisms.
\emph{For simplicity, $(AG, \rho_{AG}, [\cdot\,,\cdot]_{AG})$ will
be written $(A, \rho, [\cdot\,,\cdot])$ in the following.}

Note that if $a\in\Gamma_M(A)$, the vector field 
$a^r$ satisfies 
$a^r\sim_\tg \rho(a)\in\mx(M)$ since 
we have 
$T_g\tg a^r(g)=T_g\tg(T_{\tg(g)}r_g\,a(\tg(g)))=T_{\tg(g)}\tg\,a(\tg(g))$ for all $g\in G$.

\subsubsection{Tangent and cotangent groupoids}\label{tangentgroupoids}

Let $G$ be a Lie groupoid over $M$ with Lie algebroid $A$. 
The tangent bundle $TG$ has a natural Lie groupoid structure over $TM$. 
This structure is obtained by applying the tangent functor to each of 
the structure maps defining $G$ (source, target, multiplication, inversion and identity section). 
We refer to $TG$ with the groupoid structure over $TM$ as the \textbf{tangent groupoid} of $G$ 
\cite{Mackenzie05}. 
The set of composable pairs $(TG)_{(2)}$ of this groupoid is equal to $T(G_{(2)})$. 
For $(g,h)\in G_{(2)}$ and a 
 pair $(v_g,w_h)\in (TG)_{(2)}$, the multiplication  is 
\[v_g\star w_h:=T\mathsf m(v_g,w_h).\]

As in \cite{Mackenzie00}, we define \textbf{star vector fields on $G$}
or \textbf{star sections of $TG$} to be vector fields 
$X\in\mx(G)$ such that there exists $\bar X\in\mx(M)$ with 
$X\sim_\s\bar X$ and $\bar X\sim_\epsilon X$, i.e., 
$X$ and $\bar X$ are $\s$-related and 
$\bar X$ and $X$ are $\epsilon$-related, or $X$ restricts to $\bar X$ on $M$.
We write then $X\overset{\star}\sim_\s\bar X$.
In the same manner, we can define $\tg$-star sections,
$X\overset{\star}\sim_\tg\bar X$ with $\bar X\in\mx(M)$ and $X\in\mx(G)$.
It is easy to see that the tangent space $TG$ is spanned 
by star vector fields at each point in $G\setminus M$.
Not also that the Lie bracket of two star sections of $TG$ 
is easily seen to be a star section.

\medskip

Consider now the cotangent bundle $T^*G$. It was shown in \cite{CoDaWe87}, that $T^*G$ is a 
Lie groupoid over $A^*$. 
The source and target  of $\alpha_g\in T_g^*G$ are defined by
\[\tilde\s(\alpha_g) \in A^*_{\s(g)}, \qquad 
\tilde{\s}(\alpha_g)(a)=\alpha_g(Tl_g(a-T\tg(a))) \quad\text{ for all }\quad a\in A_{\s(g)}\]
and 
\[\tilde\tg(\alpha_g)\in A^*_{\tg(g)},\qquad  \tilde{\tg}(\alpha_g)(b)=\alpha_g(Tr_g(b)) 
\quad\text{ for all }\quad 
b\in A_{\tg(g)}.
\]
The multiplication on $T^*G$ will also be denoted by $\star$ and is defined by
\[(\alpha_g\star \beta_h)(v_g\star w_h)= \alpha_g(v_g)+ \beta_h(w_h)\]
 for $(v_g,w_h)\in T_{(g,h)}G_{(2)}$.

We refer to $T^*G$ with the groupoid structure over $A^*$ as the \textbf{cotangent groupoid} of $G$.

\subsection{The tangent  algebroid}\label{tangentalgebroids}

Let $M$ be a smooth manifold. The tangent bundle of $M$ is denoted by $p_M:TM\rmap M$. 
Consider now $q_A:A\longrightarrow M$ a vector bundle over $M$. 
The tangent bundle $TA$ has a natural structure of vector bundle over $TM$, defined 
by applying the tangent functor to each of the structure maps that define the vector 
bundle $q_A:A\to M$. This yields  a commutative diagram

\begin{align}
\begin{xy}
\xymatrix{
TA\ar[r]^{Tq_A}\ar[d]_{p_A}&TM\ar[d]^{p_M}\\
A\ar[r]_{q_A}&M
}
\end{xy}
\end{align}

\noindent where all the structure maps of $TA\to TM$ are vector bundle morphisms 
over the corresponding structure maps of $A\to M$. In the terminology of \cite{Mackenzie05,Pradines77}
 this defines a double vector bundle. 
Note finally that the zero element in the fiber of $TA$ over $v\in TM$
is $T0^Av$, i.e., the zero section is just $T0^A\in\Gamma_{TM}(TA)$.

\subsubsection{The tangent  algebroid $TA\to TM$.}
Assume now that $q_A:A\longrightarrow M$ has a Lie algebroid 
structure with anchor map $\rho:A\rmap TM$ and Lie bracket $[\cdot,\cdot]$ on $\Gamma_M(A)$. 
Then there is a Lie algebroid structure on $TA$ over $TM$ referred to as the \textbf{tangent Lie algebroid}.

In order to describe explicitly the algebroid structure on $Tq_A:TA\rmap TM$, 
we recall first that there exists a \textbf{canonical involution}
 \begin{align}\label{canonicalinvolution}
\begin{xy}
\xymatrix{
TTM\ar[r]^{J_M}\ar[d]_{p_{TM}}&TTM\ar[d]^{Tp_M}\\
TM\ar[r]_{\id_{TM}}&TM
}
\end{xy}
\end{align}
 which is given as follows \cite{Mackenzie05, Tulczyjew77}.
Elements $(\xi;v,x;m)\in TTM$, that is, 
with $p_{TM}(\xi)=v\in T_mM$ and $Tp_M(\xi)=x\in T_mM$, are considered as 
second derivatives $$\xi=\frac{\partial^2\sigma}{\partial t\partial u}(0,0),$$
where $\sigma:\R^2\to M$ is a smooth square of elements of $M$.
The notation means that $\sigma$ is first differentiated with respect to $u$, yielding a curve 
$v(t)=\frac{\partial \sigma}{\partial u}(t,0)$ in $TM$ with $\left.\frac{d}{dt}\right\an{t=0}v(t)=\xi$.
Thus, $v=\frac{\partial \sigma}{\partial u}(0,0)=p_{TM}(\xi)$
and $x=\frac{\partial \sigma}{\partial t}(0,0)=Tp_{M}(\xi)$.
The canonical involution  $J_M:TTM\to TTM$
is defined by $$J_M(\xi):=\frac{\partial^2\sigma}{\partial u\partial t}(0,0).$$
 Now we can apply the tangent functor to the anchor map $\rho:A\rmap TM$, 
and then compose with the canonical involution to obtain a bundle map $\rho_{TA}:TA\rmap TTM$ defined by
$$\rho_{TA}=J_M\circ T\rho.$$ 
 This defines the tangent anchor map. In order to define the tangent Lie bracket, we observe that every 
section $a\in\Gamma_{M}(A)$ induces two types of sections of $TA\rmap TM$. 
The first type of section is $Ta:TM\rmap TA$, called \textbf{linear} section\footnote{Note that not all linear sections
are of this type: if $\phi:TM\to A$ is a vector bundle homomorphism over $\id_M$, and $a\in\Gamma(A)$, 
then the section $a^\phi$ of $TA\to TM$ defined by 
\[a^\phi(v_m)=Ta(v_m)+\left.\frac{d}{dt}\right\an{t=0}a_m+t\cdot\phi(v_m)
\]
for all $v_m\in TM$
is also a linear section, i.e.,
$a^\phi:TM\to TA$ is a vector bundle homomorphism over $a:M\to A$.},
 which is given by applying the tangent functor to the section 
$a:M\rmap A$. The second type of section is the \textbf{core} section $\hat{a}:TM\rmap TA$, which is defined by

\begin{equation}\label{eq:coreTqA}
\hat{a}(v_m)=T_m0^A(v_m)+_{p_A}\overline{a(m)},
\end{equation}

\noindent where $0^A:M\rmap A$ denotes the zero section, and
$\overline{a(m)}=\left.\frac{d}{dt}\right\an{t=0}ta(m)\in T_{0^A_m}A$. 
As observed in \cite{MaXu98}, sections of the form $Ta$ and $\hat{a}$ generate the module 
of sections $\Gamma_{TM}(TA)$. Therefore, the tangent Lie bracket $[\cdot\,,\cdot]_{TA}$ is completely
determined by
\[\left[Ta,Tb\right]_{TA}=T[a,b], \quad \left[Ta,\hat{b}\right]_{TA}=\widehat{[a,b]}, 
\quad\left[\hat{a},\hat{b}\right]_{TA}=0\]
for all $a,b\in\Gamma(A)$,
the extension to general sections is done using the Leibniz rule with respect to the tangent anchor $\rho_{TA}$.

\bigskip

\subsubsection{The Lie algebroid of the tangent groupoid}
If $A$ is now the Lie algebroid of a Lie groupoid $G\rr M$, then we can consider 
the Lie algebroid $q_{A(TG)}:A(TG)\to TM$ of the tangent Lie groupoid $TG\rr TM$.
Since the projection $p_G:TG\to G$ is a Lie groupoid morphism, we have a Lie algebroid
morphism $A(p_G):A(TG)\to A$ over $p_M:TM\to M$ and the 
following diagram commutes.

\begin{align}
\begin{xy}
\xymatrix{
A(TG)\ar[r]^{\quad A(p_G)}\ar[d]_{q_{A(TG)}}&A\ar[d]^{q_A}\\
TM\ar[r]_{p_M}&M
}
\end{xy}
\end{align}

Let $a$ be a section of the Lie algebroid $A$, choose $v\in TM$ and consider 
the curve $\gamma:(-\varepsilon, \varepsilon)\to TG$
defined by 
\[\gamma(t)=T\Exp(ta)v\]
for $\varepsilon$ small enough.
Then we have $\gamma(0)=v$ and 
$T\s(\gamma(t))=v$ for all $t\in (-\varepsilon,\varepsilon)$. Hence, 
$\dot \gamma(0)\in A_{v}(TG)$ and we can define a \textbf{linear} section $\beta_a:TM\to A(TG)$
by 
\begin{equation}\label{beta_a}
\beta_a(v)=\left.\frac{d}{dt}\right\an{t=0}T\Exp(ta)v
\end{equation}
for all $v\in TM$. It is easy to check that $\beta_a^r\in \mx(TG)^r$
is the complete lift of $a^r$ (see \cite{Mackenzie05}). In particular the flow of $\beta_a^r$
is $TL_{Exp(\cdot a)}$, and $(\beta_a,a)$ is a morphism of vector bundles.

In the same manner, we can consider $v\in TM$, $a\in A_{p_M(v)}$
and the curve 
$\gamma:\R\to TG$
defined by 
\[\gamma(t)=v+ta,\]
where $TM$ and $A$ are seen as subsets of $TG$, $T_MG=TM\oplus A$.
We have again $\gamma(0)=v$ and $T\s(\gamma(t))=v$ for all $t$, 
which yields 
$\dot \gamma(0)\in A_v(TG)$. Given $a\in\Gamma_M(A)$, we define a \textbf{core} section $\tilde a$ of $A(TG)$ by 
\begin{equation}\label{a_tilde}
\tilde a(v)=\left.\frac{d}{dt}\right\an{t=0}v+ta(p_M(v))
\end{equation}
for all $v\in TM$.
We have for  $v_g\in T_gG$ with $T_g\tg(v_g)=v_m$:
\begin{align*}
{\tilde a}^r(v_g)&=\tilde a(v_m)\star 0_{v_g}
=\left.\frac{d}{dt}\right\an{t=0}v_g+ta^r(g).
\end{align*}

The vector bundle $A(TG)$ is spanned by the two types of sections $\beta_a$ and $\tilde a$,
for $a\in \Gamma_M(A)$, and, using the flows of 
$\beta_a^r$ and ${\tilde b}^r\in \mx^r(TG)$, it is easy to check that the equalities 
\[\left[\beta_a,\beta_b\right]_{A(TG)}=\beta_{[a,b]}, \quad 
\left[\beta_a,\tilde b\right]_{A(TG)}=\widetilde{[a,b]}, \quad 
\left[\tilde a,\tilde b\right]_{A(TG)}=0
\]
hold for all $a,b\in\Gamma_M(A)$.

 There exists a natural injective bundle map
\begin{equation}\label{i_{AG}}
\iota_{A}:A\rmap TG
\end{equation}
over $\epsilon:M\to G$.
The canonical involution $J_G:TTG\rmap TTG$ restricts to an isomorphism of Lie algebroids 
$j_G:TA\rmap A(TG)$. More precisely, there exists a commutative diagram
\begin{align}
\begin{xy}
\xymatrix{
TA\ar[r]^{j_G}\ar[d]_{T\iota_{A}}&A(TG)\ar[d]^{\iota_{A(TG)}}\\
TTG\ar[r]_{J_G}&TTG
}
\end{xy}
\end{align}
We  get easily  the following identities 
\[j_G\circ Ta=\beta_a\qquad \text{ and } \qquad j_G\circ \hat a=\tilde a,\]
where $\hat{a}$ is defined as in \eqref{eq:coreTqA}. The equality
\[\rho_{A(TG)}\circ j_G=J_M\circ T\rho_A=\rho_{TA}
\]
is  verified on linear and core sections.
This shows that the Lie algebroid $A(TG)\to TM$ of the tangent groupoid is canonically 
isomorphic to the tangent Lie algebroid $TA\to TM$ of $A$.

\bigskip

\subsubsection{The standard Lie algebroid $TA\to A$}
We study now the Lie algebroid structure $TA\to A$ in terms of special sections.
Consider a linear vector field on $A$, i.e., a section  $X$ of $TA\to A$ 
such that the map
\begin{align*}
\begin{xy}
\xymatrix{
A\ar[d]_{q_A}\ar[r]^X&TA\ar[d]^{Tq_A}\\
M\ar[r]_{\bar X}&TM
}
\end{xy}
\end{align*} 
is a vector bundle morphism over  $\bar X\in\mx(M)$. For any $b\in\Gamma_M(A)$, the core
section $b^{\uparrow}$ of $TA\to A$ associated to $b$
is defined by
\begin{equation}\label{core_uparrow}
b^{\uparrow}(a_m)=\left.\frac{d}{dt}\right\an{t=0}a_m+tb(m)
\end{equation}
for all $a_m\in A$.

The flow $\phi^X$ of $X$ is then a vector bundle 
morphism over the flow $\phi^{\bar X}$ of $\bar X$
\begin{align*}
\begin{xy}
\xymatrix{
A\ar[d]_{q_A}\ar[r]^{\phi^X_t}&A\ar[d]^{q_A}\\
M\ar[r]_{\phi^{\bar X}_t}&M
}
\end{xy}
\end{align*} 
for $t\in\R$ such that this makes sense. The flow of $b^{\uparrow}$ is given by 
$\phi^{b^{\uparrow}}_t(a_m)=a_m+tb(m)$
for all $a_m\in A$ and all $t\in\R$. In particular, for every $t\in \R$, $\phi^{b^\uparrow}_t$ is a vector bundle automorphism of $A$, covering the identity.
\begin{align*}
\begin{xy}
\xymatrix{
A\ar[d]_{q_A}\ar[r]^{\phi^{b^{\uparrow}}_t}&A\ar[d]^{q_A}\\
M\ar[r]_{\id_M}&M
}
\end{xy}
\end{align*} 

\begin{lemma}\label{linear_core_of_TA_A}
Let $A$ be a Lie algebroid. Choose linear sections $(X,\bar X), (Y,\bar Y)$ of
$TA\to A$ and sections $a,b\in\Gamma_M(A)$. Then 
\[\left([X,Y],\left[\bar X,\bar Y\right]\right)\]
is a linear section of $TA\to A$,
\[\left[a^\uparrow, b^\uparrow\right]=0\]
and \[\left[X, a^\uparrow\right]=(D_Xa)^\uparrow,\]
where 
\[D_Xa\in\Gamma(A), \quad (D_Xa)(m)=\left.\frac{d}{dt}\right\an{t=0}
\underset{\in A_m}{\underbrace{\left(\phi^X_{-t}\circ a\circ \phi^{\bar X}_t\right)(m)}}\]
for all $m\in M$.
\end{lemma}
\begin{proof}
This is easy to show using the flows of linear and core sections, see also \cite{Mackenzie92}. 
\end{proof}

\bigskip

Recall that if $G\rr M$ is a Lie groupoid, then  the tangent space $TG$ is spanned outside of $M$
by star sections $X\overset{\star}\sim_\s\bar X$, i.e. with $X\in\mx(G)$, $\bar X\in\mx(M)$ such that 
$X\sim_\s\bar X$ and $\bar X\sim_\epsilon X$ (see Subsection \ref{tangentgroupoids}).
The bracket of two star sections is a star section and 
for any $a_m\in A$, the properties of 
the star vector field  $X\in\mx(G)$ over $\bar X\in\mx(M)$ 
imply that $(TX)(a_m)$ has value in $A_{\bar X(m)}(TG)$.
Hence, we can  mimic the construction 
of the Lie algebroid map associated to a Lie groupoid morphism and  we can consider the 
map $A(X):A\to A(TG)$, $A(X)(a_m)=T_mX(a_m)$ over $\bar X:M\to TM$.
In the same manner, if $a\in\Gamma_M(A)$, then 
we can define $\bar a: A\to A(TG)$ by 
$\bar a(b_m)=\beta_b(0^{TM}_m)+_{q_{A(TG)}}\tilde a(0^{TM}_m)$
for any section $b\in\Gamma_M(A)$ such that $b(m)=b_m$.
We have $j_G\inv\circ \bar a=a^\uparrow:A\to TA$
and for a star section $X\overset{\star}\sim_\s\bar X$ of $TG\to TM$, the map
$\tilde X:=j_G\inv\circ A(X):A\to TA$ 
 is a linear section $(\tilde X,\bar X)$ of $TA\to A$. There is a unique Lie algebroid 
structure on $A(TG)$ over $A$ making $j_G:TA\to A(TG)$ into a Lie algebroid isomorphism over the identity on $A$ \cite{Mackenzie92, MaXu98, Mackenzie00}.

\subsection{Flat partial connections}
The following definition will be crucial in this paper.
\begin{definition} (\cite{Bott72})
Let $M$ be a smooth manifold and $F\subseteq TM$ a smooth involutive vector subbundle
of the tangent bundle.
 Let $E\to M$ be a vector bundle over $M$.
A \textbf{$F$-partial connection} 
is a map $\nabla:\Gamma( F)\times \Gamma(E)\to\Gamma(E)$, written 
$\nabla(X,e)=:\nabla_Xe$ for $X\in\Gamma(F)$ and $e\in\Gamma(E)$,
such that:
\begin{enumerate}
\item $\nabla$ is tensorial in the $F$-argument,
\item $\nabla$ is $\R$-linear in the $E$-argument,
\item $\nabla$ satisfies the Leibniz rule
\[\nabla_{X}(fe)=X(f)e+f\nabla_{X}e\]
for all $X\in\Gamma(F)$, $e\in\Gamma(E)$, $f\in C^\infty(M)$.
\end{enumerate}
The connection is \textbf{flat} if 
\[\nabla_{[X,Y]}e=\left(\nabla_X\nabla_Y-\nabla_Y\nabla_X\right)e
\]
for all $X,Y\in\Gamma(F)$ and $e\in\Gamma(E)$.
\end{definition}

\begin{example}[The Bott connection]
Let $M$ be a smooth manifold and $F\subseteq TM$
an involutive subbundle.
The \textbf{Bott connection} 
\[\nabla^F:\Gamma(F)\times\Gamma(TM/F)\to\Gamma(TM/F)\]
defined by 
\begin{equation*}
\nabla_{X}\bar Y=\overline{[X,Y]},
\end{equation*}
where $\bar Y\in\Gamma(TM/F)$ is the projection 
of $Y\in\mx(M)$,
is a flat $F$-partial connection on $TM/F\to M$.

\medskip

The class $\bar Y\in\Gamma(TM/F)$ of a vector field 
is $\nabla^F$-parallel if and only if $[Y,\Gamma(F)]\subseteq \Gamma(F)$. Since $F$ is involutive, this does not
depend on the representative of $\bar Y$. We say
by abuse of notation that 
$Y$ is $\nabla^F$-parallel.
\end{example}

The following proposition can be easily shown by using the fact that the parallel transport
defined by a flat connection does not depend 
on the chosen path in simply connected sets
(see \cite{Jotz11b}, \cite{JoRaZa11} for similar statements).
\begin{proposition}\label{parallel_sections}
Let $E\to M$ be a smooth vector bundle of rank $k$,  $F\subseteq TM$
an involutive subbundle and $\nabla$ a flat partial $F$-connection on $E$.
Then there exists for each point $m\in M$ 
 a frame of local $\nabla$-parallel sections $e_1,\ldots,e_{k}
\in\Gamma(E)$ defined on an open neighborhood $U$ of $m$ in $M$.
\end{proposition}

We have also the following lemma. The first statement 
is a straightforward consequence of the Leibniz identity
and the proof of the second statement can be checked easily 
in coordinates adapted to the foliation.
\begin{lemma}\label{immediate_lemma}
Let $E\to M$ be a smooth vector bundle of rank $k$,  $F\subseteq TM$
an involutive subbundle and $\nabla$ a partial $F$-connection on $E$.
\begin{enumerate}
\item Assume that $f\in C^\infty(M)$ is $F$-invariant, i.e.,
$X(f)=0$ for all $ X\in \Gamma(F)$. Then 
$f\cdot e$ is $\nabla$-parallel for any 
$\nabla$-parallel section $e\in\Gamma(E)$. 
\item  Assume that the foliation defined by $F$ on $M$ is simple, i.e.
the leaf space has a smooth manifold structure such that the quotient 
$\pi:M\to M/F$ is a smooth surjective submersion. Then $X\in\mx(M)$ is $\nabla^F$-parallel
if and only if there exists $\bar X\in\mx(M/F)$ such that $X\sim_\pi\bar X$.
\end{enumerate}
\end{lemma}

\vspace*{0.5cm}

\section{Foliated groupoids}\label{fol_gpds}
\subsection{Definition and properties}
\begin{definition}\label{def_fol_gpd}
Let $G\rr M$ be a Lie groupoid. A subbundle $F_G\subseteq TG$ is 
\textbf{multiplicative} if it is a subgroupoid  of $TG\rr TM$ over
$F_G\cap TM=:F_M$. We say that $F_G$ is a \textbf{multiplicative foliation} on $G$
if it is involutive and multiplicative.
The pair $(G\rr M, F_G)$ is then called a \textbf{foliated groupoid}.
\end{definition}

\begin{remark}
Multiplicative subbundles were introduced in \cite{Tang06}
as follows. A subbundle $F_G\subseteq TG$ is multiplicative
if for all composable $g,h\in G$ and $u\in F_G(g\star h)$, there exist
$v\in F_G(g)$, $w\in F_G(h)$
such that $u=v\star w$.
It is easy to check that a multiplicative foliation in the sense of Definition 
\ref{def_fol_gpd} is multiplicative in the sense of \cite{Tang06}, but the converse is not necessarily true, 
unless for instance if the Lie groupoid is a Lie group (see \cite{Jotz11a}). 
The case of involutive wide subgroupoids of $TG\rr TM$ has also been studied in  \cite{Behrend05}.
\end{remark}

The following result about multiplicative subbundles can be found in \cite{Jotz11b}.

\begin{lemma}\label{constant_rank}
Let $G\rr M$ be a Lie groupoid and
 $F_G\subseteq TG$ a multiplicative subbundle. 
 Then the intersection $F_M:=F_G\cap TM$ has constant rank on $M$.
Since it is the set of units of $F_G$ seen as a
subgroupoid of $TG$, 
the pair $F_G\rr F_M$
is a Lie groupoid. 

The bundle $F_G\an{M}$
splits as $F_G\an{M}=F_M\oplus F_\core$, where 
$F_\core:=F_G\cap A$.
We have \[(F_G\cap T^\s G)(g)= F_\core(\tg(g))\star 0_g
=T_{\tg(g)}r_g\left(F_\core(\tg(g))\right)\] for all
$g\in G$.

In the same manner, if
$F^\tg:=(F_G\cap T^\tg G)\an{M}$, we have 
$(F_G\cap T^\tg G)(g)=0_g\star F^\tg(\s(g))$ for all
$g\in G$.

As a consequence, the intersections $F_G\cap T^\tg G$ and 
$F_G\cap T^\s G$ have  constant rank on $G$.
\end{lemma}

The previous lemma says that a multiplicative subbundle $F_G\subseteq TG$ determines a $\mathcal{VB}$-groupoid 
\begin{align*}
\begin{xy}
\xymatrix{
F_G\ar[r]^{p_G}\ar@<.6ex>^{T\s}[d]\ar@<-.6ex>_{T\tg}[d]& G\ar@<.6ex>^{\s}[d]\ar@<-.6ex>_{\tg}[d] \\
F_M\ar[r]_{p_M}&M}
\end{xy}
\end{align*} 
with  \textbf{core} $F_\core$.
As a corollary, one gets  that the source and target maps  of $F_G\rr F_M$ are fiberwise surjective.
\begin{corollary}\label{surjectivity}
Let $G\rr M$ be a  Lie groupoid 
and $F_G$ a multiplicative subbundle of $TG$.
 The induced
maps $T_g\s:F_G(g)\to F_G(\s(g))\cap T_{\s(g)}M
$ and $T_g\tg:F_G(g)\to F_G(\tg(g))\cap T_{\tg(g)}M$
are surjective for each $g\in G$.
\end{corollary}

\begin{proposition}\label{subalgebroids}
Let $F_G$ be a multiplicative foliation on a Lie groupoid $G\rr M$. 
Then $F_\core$ is a subalgebroid of $A$  and $F_M$ is an involutive  subbundle of $TM$.
\end{proposition}
\begin{proof}
Choose first two sections $a,b\in\Gamma(F_\core)$. Then we have 
$a^r, b^r\in\Gamma(F_G)^r$ by Lemma \ref{constant_rank}
and hence $[a,b]^r=\left[a^r,b^r\right]\in\Gamma(F_G)$
since $F_G$ is involutive. Again by Lemma \ref{constant_rank}, we find 
thus that $[a,b]\in\Gamma(F_\core)$.

\medskip

Choose now two sections $\bar X,\bar Y\in\Gamma(F_M)$ and 
$X,Y\in\Gamma(F_G)$
that are $\s$-related to $\bar X,\bar Y$. This is possible by Corollary
\ref{surjectivity}.
We find then 
that $[X,Y]\in\Gamma(F_G)$ and $[X,Y]\sim_\s[\bar X,\bar Y]$.
Hence, $[\bar X,\bar Y]\in\Gamma(F_M)$.
\end{proof}


As explained in \cite{Jotz11c, Ortiz08t}, a multiplicative foliation is a $\mathcal LA$-groupoid
in the sense of \cite{Mackenzie00},
\begin{displaymath}\begin{xy}
\xymatrix{ 
 F_G\ar@<.6ex>^{T\s}[rr]\ar@<-.6ex>_{T\tg}[rr]
\ar[dd]_{p_G}\ar[dr]^{\iota_{F_G}}&&F_M\ar@{-->}[dd]\ar[dr]^{\iota_{F_M}}&\\
& TG \ar@<.6ex>^{\qquad T\s}[rr]\ar@<-.6ex>_{\qquad T\tg}[rr]\ar[dl]_{p_G}&&TM\ar[dl]^{p_M}\\
G\ar@<.6ex>^{\s}[rr]\ar@<-.6ex>_{\tg}[rr]&&M&
}
\end{xy}
\end{displaymath}
where and $\iota_{F_M}, \iota_{F_G}$ are the inclusions 
of $F_M$ in $TM$ and of $F_G$ in $TG$, respectively.
That is, 
$(F_M\to M, \iota_{F_M}, [\cdot\,,\cdot])$ is a Lie algebroid over $M$
and the quadruple $(F_G; G, F_M; M)$ is such that 
$F_G$ has  both a  Lie groupoid structure over $F_M$ and a Lie algebroid
structure   over $G$ such that the two structures on $F_G$ are compatible in the sense 
that the maps defining  the groupoid structure are all Lie algebroid 
morphisms. 
Furthermore, the double source map 
\[(p_G, T\s):F_G\to G\times_M F_M=\{(g,v)\in G\times F_M\mid v\in F_M(\s(g))\} 
\] 
is a surjective submersion by Corollary \ref{surjectivity}.

\subsection{The connection associated to a foliated groupoid}

 Recall that 
if $G\rr M$ is a Lie groupoid and $F_G$ is a multiplicative foliation on $G$,
then $F_M:=F_G\cap TM$ and $F_\core=F_G\cap A$ 
are subalgebroids of $TM\to M$ and $A\to M$, respectively (Proposition \ref{subalgebroids}).

In the main theorem of this subsection, we show that 
if $F_G$ is a multiplicative foliation on the Lie groupoid $G\rr M$, then the Bott $F_G$-connection 
on $TG/F_G$ induces  a well-defined partial $F_M$-connection  on $A/F_\core$.
We write $\bar a$ for the class in $A/F_\core$ of $a\in\Gamma(A)$.
\begin{theorem}\label{the_connection}
Let $(G\rr M,F_G)$ be a  Lie groupoid endowed with a multiplicative foliation.
Then there is a partial $F_M$-connection on $ A/F_\core$ 
\begin{align}
\nabla:\Gamma(F_M)\times\Gamma\left(A/F_\core\right)&\to \Gamma\left(A/F_\core\right).
\end{align}
with the following properties:
\begin{enumerate}
\item $\nabla$ is flat.
\item If $a\in\Gamma(A)$ is $\nabla$-parallel, i.e., 
$\nabla_{\bar X}\bar a=0$ for all $\bar X\in\Gamma(F_M)$, then $[a,b]\in\Gamma(F_\core)$ 
for all $b\in\Gamma(F_\core)$.
\item If $a,b\in\Gamma(A)$ are $\nabla$-parallel, then $[a,b]$ is also $\nabla$-parallel.
\item If $a\in\Gamma(A)$ is  $\nabla$-parallel, then $[\rho(a),\bar X]\in\Gamma(F_M)$ for all
$\bar X\in\Gamma(F_M)$. That is, $\rho(a)\in\mx(M)$ is $\nabla^{F_M}$-parallel.
\end{enumerate}
Furthermore,  $a\in\Gamma(A)$ is  $\nabla$-parallel if and only if the right invariant vector 
field $a^r\in\mx(G)$ is parallel with respect to the Bott connection defined by $F_G$, that is
\[\left[a^r,\Gamma(F_G)\right]\subseteq \Gamma(F_G).\]
\end{theorem}

In the following, we call 
a vector field $X\in\mx(G)$ a
\textbf{$\tg$-section} if 
there exists $\bar X\in\mx(M)$
such that $X\sim_\tg\bar X$.
Similarly, a one-form
$\eta\in\Omega^1(G)$ is a $\tg$-section of $T^*G$ if 
$\tilde\tg\circ\eta=\bar\eta\circ\tg$
for some $\bar\eta\in\Gamma(A^*)$.
Analogously, we can define 
\textbf{$\s$-sections} of $TG$ and $T^*G$.

It is easy to see that
$F_G$ and 
its annihilator ${F_G}^\circ\subseteq T^*G$ are spanned by their $\tg$-sections.

\begin{lemma}
Let $(G\rr M,F_G)$ be a Lie groupoid endowed with a multiplicative subbundle $F_G\subseteq TG$.
Choose $\bar X\in \Gamma(F_M)$, $\bar\eta\in\Gamma(F_G^\circ\cap A^*)$ and 
$\tg$-sections $X\sim_\tg \bar X$,
$\eta\sim_{\tilde\tg}\bar\eta$ of $F_G$ and $F_G^\circ$, respectively.
Then the identity
\begin{equation}\label{key}
\eta(\ldr{a^r}X) =\tg^*\epsilon^*\bigl(
\bar\eta(\ldr{a^r}X)\bigr)
\end{equation}
holds for any section $a\in\Gamma(A)$.
\end{lemma}

\begin{proof}
Choose $g\in G$ and set $p=\tg(g)$. 
For all $t\in (-\varepsilon, \varepsilon)$ for a small $\varepsilon$, 
we have 
\begin{align*}
X(\Exp(ta)(p)\star g)&=
(X(\Exp(ta)(p)))\star \bigl(X(\Exp(ta)(p))\bigr)\inv
\star X(\Exp(ta)(p)\star g).
\end{align*}
The vector  $$\bigl(X(\Exp(ta)(p))\bigr)\inv
\star X(\Exp(ta)(p)\star g)$$
 is an element 
of $F_G(g)$ for all $t\in (-\varepsilon, \varepsilon)$
and will be written $v_t(g)$ to simplify the notation.
We compute 
\begin{align*}
\eta(\ldr{a^r}X)(g)
&= \eta(g)
\left(\left.\frac{d}{dt}\right\an{t=0}\bigl(T_{L_{\Exp(ta)}(g)}L_{\Exp(-ta)}X(\Exp(ta)(p)\star
g)\right)\\
&=\left.\frac{d}{dt}\right\an{t=0}
 (\bar\eta(p)\star \eta(g)) \Bigl(\bigl(T_{\Exp(ta)}L_{\Exp(-ta)}X(\Exp(ta)(p))\bigr)\star v_t(g)\Bigr)\\
&=\left.\frac{d}{dt}\right\an{t=0}
 \bar\eta(p)\bigl(T_{\Exp(ta)}L_{\Exp(-ta)}X(\Exp(ta)(p))\bigr)
+\left.\frac{d}{dt}\right\an{t=0}\eta(g)(v_t(g))\\
&=\bar\eta(\ldr{a^r}X)(p)+
\left(\left.\frac{d}{dt}\right\an{t=0}0\right)
=\bar\eta(\ldr{a^r}X)(\tg(g)).
\end{align*}
\end{proof}

Let $G\rr M$ be a  Lie groupoid and $F_G$ a multiplicative subbundle of $TG$.
Recall that
the restriction to $M$ of $F_G$ splits as 
\begin{equation}\label{splitting_of_F}
F_G\an{M}=F_M\oplus F_\core,
\end{equation}
where 
$F_\core=F_G\cap A$ and 
$F_M=F_G\cap TM$. Recall also the notation $F^\tg=F_G\cap T_M^\tg G$.

\begin{theorem}\label{lie_der_of_xi_section}
Let $(G\rr M, F_G)$ be a Lie  groupoid endowed with a multiplicative subbundle $F_G\subseteq TG$,
$X$ a $\tg$-section of
$F_G$, i.e., $\tg$-related to some
$\bar X\in\Gamma(F_M)$,
and consider $a\in\Gamma(A)$. Then the derivative 
$\ldr{a^r}X$ 
can be written 
as a sum 
\[\ldr{a^r}X=Z_{a,X}+ b_{a,X}^r\]
with 
$b_{a,X}\in\Gamma(A)$,
and $Z_{a,X}$ a  $\tg$-section of $F_G$. 
In addition, if $X\sim_\tg0$, 
then $\ldr{a^r}X\in\Gamma(F_G\cap T^\tg G)$. In particular, its restriction to $M$ is a section of $F^\tg$
and $b_{a,X}$ is a section of $F_\core$.
\end{theorem}

\begin{proof}
Set \[b_{a,X}(m)=\ldr{a^r}X(m)-T_m\s\left(\ldr{a^r}X(m)\right)\]
for all $m\in M$ and
\[Z_{a,X}:= \ldr{a^r}X-b_{a,X}^r.\]
In particular, $Z_{a,X}(m)=T_m\s\left(\ldr{a^r}X(m)\right)$ for all $m\in M$. First, we see 
that $Z_{a,X}$ is a $\tg$-section of $TG$ since
\[Z_{a,X}=\ldr{a^r}X-b_{a,X}^r\sim_\tg \left[\rho(a),\bar X\right]-\rho(b_{a,X}).\]

Recall that ${F_G}^\circ\subseteq T^*G$  is spanned by its $\tg$-sections. Hence,
to show that $Z_{a,X}\in\Gamma(F_G)$,
 it suffices 
to show that $\eta(Z_{a,X})=0$ for all 
$\tg$-sections $\eta$ of ${F_G}^\circ$.

Setting $\tg(g)=p$, we have immediately
\begin{align*}
\eta(Z_{a,X})(g)&\,\,=\,\,\eta(\ldr{a^r}X)(g)-\bar\eta(b_{a,X})(p)\\
&\overset{\eqref{key}}=\bar\eta(\ldr{a^r}X)(p)-\bar\eta(b_{a,X})(p)\\
&\,\, =\,\,\bar\eta(T_p\s(\ldr{a^r}X(p)))
=0
\end{align*}
since $\bar\eta$ vanishes on $TM$. This shows that $Z_{a,X}$ is a section of $F_G$.
\medskip

Assume now that $X\sim_\tg 0$. This means that $X\in\Gamma(F_G\cap T^\tg G)$
and $X$ can be written
$X=\sum_{i=1}^rf_i\cdot b_i^l$
with $f_1,\ldots,f_r\in C^\infty(G)$ and 
$b_1,\ldots,b_r\in\Gamma(F^\tg)$.
We get then 
\[\ldr{a^r}X=\ldr{a^r}\left(\sum_{i=1}^rf_i\cdot b_i^l\right)
=\sum_{i=1}^ra^r(f_i)b_i^l+\sum_{i=1}^rf_i\cdot\left[a^r, b_i^l\right]
=\sum_{i=1}^ra^r(f_i)b_i^l,
\] 
which is again a section of $F_G\cap T^\tg G$.
\end{proof}

\bigskip

Assume now that $F_G$ is involutive and define \[\nabla:\Gamma(F_M)\times\Gamma(A/F_\core)\to \Gamma(A/F_\core)\]
by 
 \begin{equation*}
\nabla_{\bar X}\bar a=-\overline{b_{a,X}},
\end{equation*}
with $b_{a,X}$ as in Theorem \ref{lie_der_of_xi_section},
for any choice of $\tg$-section $X\in\Gamma(F_G)$ such that
$X\sim_\tg\bar X$ and any choice of representative $a\in\Gamma(A)$ for $\bar a$.
We will show that this is a well-defined partial 
$F_M$-connection and complete the proof of Theorem \ref{the_connection}.
\begin{proof}[Proof of Theorem \ref{the_connection}]
Choose $X,X'\in\Gamma(F_G)$ such that $X\sim_\tg\bar X$ and 
$X'\sim_\tg \bar X$. Then $Y:=X-X'\sim_\tg 0$
and, by Theorem \ref{lie_der_of_xi_section}, we find $b_{a,Y}\in\Gamma(F_\core)$
for any $a\in\Gamma(A)$, i.e., $\overline{b_{a,X}}=\overline{b_{a,X'}}$.

Choose now $a\in\Gamma(F_\core)$ and $X\in\Gamma(F_G)$, $X\sim_\tg\bar X\in\Gamma(F_M)$. 
Then we have $a^r\in\Gamma(F_G)$ and since $F_G$ is involutive,
$\ldr{a^r}X\in\Gamma(F_G)$. Again, since $Z_{a,X}\in \Gamma(F_G)$, we find $b_{a,X}\in\Gamma(F_\core)$.
This shows that $\nabla$ is well-defined.

\medskip

By definition, 
if $a\in\Gamma(A)$ is such that $\nabla_{\bar X}\bar a=0$ for all
$\bar X\in\Gamma(F_M)$, then we have 
$\ldr{a^r}X=Z_{a,X}+b_{a,X}^r\in \Gamma(F_G)$
for all $\tg$-descending sections $X\in\Gamma(F_G)$. 
Since $\Gamma(F_G)$ is spanned as a $C^\infty(G)$-module
by its $\tg$-descending sections, we get
\[\left[a^r,\Gamma(F_G)\right]\subseteq \Gamma(F_G).\]
Conversely, $\left[a^r,\Gamma(F_G)\right]\subseteq \Gamma(F_G)$ implies immediately
$ \nabla_{\bar X}\bar a=0$ for all
$\bar X\in\Gamma(F_M)$.
This proves the second claim of the theorem.

\medskip

We check that $\nabla$ is a flat partial $F_M$-connection. 
Choose $a\in\Gamma(A)$, $\bar X\in\Gamma(F_M)$, $X\in\Gamma(F_G)$ such that $X\sim_\tg\bar X$
and $f\in C^\infty(M)$. Then we have  $\tg^*f\cdot X\sim_\tg f\bar X$ and 
\begin{align*}
\ldr{a^r}(\tg^*f\cdot X)=\tg^*(\rho(a)(f))\cdot X+ \tg^*f\cdot \ldr{a^r}X.
\end{align*}
In particular, we find 
\begin{align*}
b_{a,\tg^*f\cdot X}&=(1-T\s)\left(\tg^*(\rho(a)(f))\cdot X+\tg^*f\cdot \ldr{a^r} X\right)\an{M}\\
&=\rho(a)(f)\cdot (1-T\s)X\an{M}+f\cdot (1-T\s)(\ldr{a^r}X)\an{M}.
\end{align*}
Since $(T\s-1)X\an{M}\in\Gamma(F_\core)$, this leads to 
$\overline{b_{a,\tg^*f \cdot X}}=f\cdot \overline{b_{a,X}}$ and hence 
$\nabla_{f\bar X}\bar a =-f\cdot  \overline{b_{a,X}}
=f\cdot\nabla_{\bar X}\bar a$.

Since $(fa)^r=\tg^* f\cdot a^r$, we have in the same manner
\begin{align*}
\ldr{(fa)^r}X=-\ldr{X}(\tg^* f\cdot a^r)=-\tg^*(\bar X(f))\cdot a^r+\tg^*f\cdot \ldr{a^r}X,
\end{align*}
which leads to $\nabla_{\bar X}(f\cdot \bar a)=\bar X(f)\cdot \bar a+f\cdot \nabla_{\bar X}\bar a$.

Choose $\bar X,\bar Y\in\Gamma(F_M)$ and $X,Y\in \Gamma(F_G)$ such that 
$X\sim_\tg\bar X$ and $Y\sim_\tg\bar Y$. Then we have 
$[X,Y]\sim_\tg[\bar X,\bar Y]$ and $[X,Y]\in\Gamma(F_G)$ since $F_G$ is involutive.
For any $a\in\Gamma(A)$, we have by the Jacobi-identity:
\begin{align*}
\ldr{a^r}[X,Y]&=\left[\ldr{a^r}X,Y\right]-\left[\ldr{a^r}Y,X\right]\\
&=\left[Z_{a,X} + b_{a,X}^r,Y\right]-\left[Z_{a,Y} + b_{a,Y}^r,X\right]\\
&=\left[Z_{a,X},Y\right]-\left[Z_{a,Y},X\right] +\ldr{b_{a,X}^r}Y-\ldr{b_{a,Y}^r}X\\
&=\left[Z_{a,X},Y\right]-\left[Z_{a,Y},X\right] 
+Z_{b_{a,X},Y}+{b_{b_{a,X},Y}}^r-Z_{b_{a,Y},X}-{b_{b_{a,Y},X}}^r.
\end{align*}
Since $\left[Z_{a,X},Y\right]-\left[Z_{a,Y},X\right]+Z_{b_{a,X},Y}-Z_{b_{a,Y},X}$
is a $\tg$-section of $F_G$, we find that
\[\nabla_{[\bar X,\bar Y]}\bar a=
\overline{b_{b_{a,Y},X}}-\overline{b_{b_{a,X},Y}}
=\nabla_{\bar X}\nabla_{\bar Y}\bar a-\nabla_{\bar Y}\nabla_{\bar X}\bar a\]
which shows the flatness of $\nabla$.

Choose now $a\in\Gamma(A)$ such that 
$\nabla_{\bar X}\bar a=0\in\Gamma(A/F_\core)$ for all $\bar X\in\Gamma(F_M)$.
If $b\in\Gamma(F_\core)$, then $b^r\in\Gamma(F_G)$, $\rho(b)\in\Gamma(F_M)$ and 
$b^r\sim_\tg\rho(b)$. This leads to
$$\overline{[b,a]}=\nabla_{\rho(b)}\bar a=0\in\Gamma(A/F_\core)$$
and hence $[a,b]\in\Gamma(F_\core)$. This shows $2$.
For each $\bar X\in\Gamma(F_M)$, there exists 
$X\in\Gamma(F_G)$ such that $X\sim_\tg\bar X$.
Since $[a^r,X]\in\Gamma(F_G)$,  $a^r\sim_\tg\rho(a)$ and $T\tg(F_G)=F_M$, we find 
$[\rho(a),\bar X]\in\Gamma(F_M)$, which proves $4$.

To show $3.$, choose two sections $a,b\in\Gamma(A)$ such that
$\bar a$ and $\bar b$ are $\nabla$-parallel.
We have then for any $\tg$-section $X\sim_\tg\bar X$ of $F_G$: 
\begin{align*}
\ldr{[a,b]^r}X&=\ldr{a^r}(Z_{b,X}+b_{b,X}^r)-\ldr{b^r}(Z_{a,X}+b_{a,X}^r)\\
&=\ldr{a^r}(Z_{b,X})+[a,b_{b,X}]^r-\ldr{b^r}(Z_{a,X})-[b,b_{a,X}]^r.
\end{align*}
Since $\bar a$ and $\bar b$ are $\nabla$-parallel, this yields 
$\nabla_{\bar X}\overline{[a,b]}=-\overline{[a,b_{b,X}]}+\overline{[b,b_{a,X}]}$.
Since $\bar a$ and $\bar b$ are parallel, 
we have 
$b_{b,X}, b_{a,X}\in\Gamma(F_\core)$ and $3.$ follows using $2$.
\end{proof}

\begin{remark}
Let $(G\rr M, F_G)$ be a foliated Lie groupoid and consider the multiplicative 
Dirac structure $F_G\oplus {F_G}^\circ$ on $G$. The ``tangent'' part of the Courant 
algebroid $\lie B( F_G\oplus {F_G}^\circ)$ in \cite{Jotz11c} is in this case just
$(F_M\oplus A)/F^\tg$, which is isomorphic as a vector bundle 
to $F_M\oplus(A/F_\core)$. The induced bracket
on sections of $F_M\oplus(A/F_\core)$
defines a map $\Gamma(F_M)\times\Gamma(A/F_\core)\to \Gamma(A/F_\core)$, 
\[\left(\bar X,\bar a\right)\mapsto \pr_{A/F_\core}\left[\left(\bar X,0\right), \left(0,\bar a\right)\right],
\]
which can be checked to be  exactly the connection $\nabla$. 
\end{remark}

\subsection{Involutivity of a multiplicative subbundle of $TG$}
It is natural to ask here how the involutivity of $F_G$ is encoded in
the data $(F_M,F_\core,\nabla)$.
For an arbitrary (not necessarily involutive) multiplicative subbundle $F_G\subseteq TG$,
we can consider the map
\[\tilde\nabla:\Gamma(F_M)\times\Gamma(A)\to \Gamma(A/F_\core),\] 
 \begin{equation*}
\tilde\nabla_{\bar X}a=-\overline{b_{a,X}}
\end{equation*}
which is well-defined by the proof of Theorem \ref{the_connection}.

\begin{theorem}\label{inv_crit}
Let $(G,F_G)$ be a source-connected Lie groupoid endowed with a multiplicative subbundle.
Then $F_G$ is involutive if and only if the following holds:
\begin{enumerate}
\item $F_M\subseteq TM$  is involutive,
\item $\tilde\nabla$ vanishes on sections of $F_\core$,
\item the induced map $\nabla:\Gamma(F_M)\times\Gamma(A/F_\core)\to \Gamma(A/F_\core)$ is a flat partial
 $F_M$-connection on $A/F_\core$.
\end{enumerate}
\end{theorem}

The proof of this theorem is a simplified version of the proof of the general
criterion for the closedness of multiplicative Dirac structures (see \cite{thesis}).
\begin{proof}
We have already shown in Proposition \ref{subalgebroids} and  Theorem \ref{the_connection}
that the involutivity of $F_G$ implies (1), (2) and (3).

For the converse implication, 
recall that the  $\tg$-star sections of $F_G$ span $F_G$ outside 
of the set of units $M$. Hence, it is sufficient 
to show involutivity on $\tg$-star sections and right-invariant sections of $F_G$.
Choose first two right-invariant sections $a^r,b^r$ of $F_G$, i.e. with $a,b\in\Gamma(F_\core)$.
Then we have $\rho(b)\in\Gamma(F_M)$, $b^r\sim_\tg\rho(b)$  and, since 
$\tilde\nabla_{\rho(b)}a=0$ by (2), we find that $[a^r,b^r]\in\Gamma(F_G)$. 
In the same manner, by the definition of $\tilde\nabla$ and Theorem \ref{lie_der_of_xi_section},
Condition (2) implies that the bracket of a right-invariant section of $F_G$ and a $\tg$-star section is always
a section of $F_G$.

We have thus only to show that the bracket of two $\tg$-star sections of $F_G$ is again a section of $F_G$.
Let $F_\core^\circ$ be the annihilator of $F_\core$ in $A^*$ and  consider the dual $F_M$-connection on
$\left(A/F_\core\right)^*\simeq{F_\core}^\circ\subseteq A^*$, i.e. the (by (3))
flat connection 
\[\nabla^*:\Gamma(F_M)\times \Gamma({F_\core}^\circ)\to \Gamma({F_\core}^\circ)
\]
given 
by 
\[\left(\nabla^*_{\bar X}\alpha\right)(\bar a)=\bar X(\alpha(\bar a))-\alpha\left(\nabla_{\bar X}\bar a\right)\]
for all $\bar X\in\Gamma(F_M)$, $\alpha\in \Gamma({F_\core}^\circ)$ and $ a\in\Gamma(A)$.

Choose  a $\tg$-star section $X\in\Gamma(F_G)$, $X\sim_\tg \bar X$,
$X\an{M}=\bar X$ and a $\tilde\tg$-section $\eta\in\Gamma(F_G^\circ)$, $\eta\sim_{\tilde\tg}\bar\eta$. 
Then, for any section $a$ of $A$, we have 
\begin{align*}
(\ldr{X}\eta)(a^r)&=X\left(\eta(a^r)\right)+\eta\left(\ldr{a^r}X\right)\\
&=\tg^*\left(\bar X(\bar \eta(\bar a))\right)+\eta\left(Z_{a,X}+ b_{a,X}^r
\right)\qquad\text{ by Theorem  \ref{lie_der_of_xi_section}}\\
&=\tg^*\left(\bar X(\bar \eta(\bar a))-\bar\eta\left(\nabla_{\bar X}\bar a\right)\right)        
\qquad\text{ since } \eta\in\Gamma(F_G^\circ) \text{ and }  Z_{a,X}\in\Gamma(F_G) \\
&=\tg^*(\nabla^*_{\bar X}\bar\eta(\bar a)).
\end{align*}
This shows that $\ldr{X}\eta\sim_{\tilde\tg}\nabla_{\bar X}^*\bar\eta\in\Gamma(F_\core^\circ)$.
Note that we have not shown yet that $\ldr{X}\eta$ is a section of $F_G^\circ$.
Choose a second $\tg$-star section $Y$ of $F_G$, $Y\sim_\tg\bar Y$ and $Y\an{M}=\bar Y$. An easy computation 
using $\eta(X)=\eta(Y)=0$ yields
\[-2\dr\left(\eta([X,Y])\right)=\ldr{X}\ldr{Y}\eta-\ldr{Y}\ldr{X}\eta-\ldr{[X,Y]}\eta.\]
Hence, we get for $a\in\Gamma(A)$:
\begin{align*}
-2\cdot a^r (\eta([X,Y]))=&\left(\ldr{X}\ldr{Y}\eta-\ldr{Y}\ldr{X}\eta-\ldr{[X,Y]}\eta\right)(a^r)\\
=&X(\ldr{Y}\eta(a^r))-\ldr{Y}\eta([X, a^r])
-Y(\ldr{X}\eta(a^r))-\ldr{X}\eta([Y, a^r])\\
&-[X,Y](\eta(a^r))+\eta([[X,Y], a^r])\\
=&\tg^*\bar X(\nabla_{\bar Y}^*\bar\eta(\bar a))+(\ldr{Y}\eta)(Z_{a,X}+b_{a,X}^r)\\
&-\tg^*\bar Y(\nabla_{\bar X}^*\bar\eta(\bar a))-(\ldr{X}\eta)(Z_{a,Y}+b_{a,Y}^r)\\
&-\tg^*[\bar X,\bar Y](\bar \eta(\bar a))-\eta\left([\ldr{a^r}X,Y]+[X,\ldr{a^r}Y]\right)\\
=&\tg^*\bar X(\nabla_{\bar Y}^*\bar\eta(\bar a))+\ldr{Y}\eta(Z_{a,X})-\tg^*(\nabla^*_{\bar Y}\bar\eta)(\nabla_{\bar X}\bar a)\\
&-\tg^*\bar Y(\nabla_{\bar X}^*\bar\eta(\bar a))-\ldr{X}\eta(Z_{a,Y})+\tg^*(\nabla^*_{\bar X}\bar\eta)(\nabla_{\bar Y}\bar a)\\
&-\tg^*[\bar X,\bar Y](\bar \eta(\bar a))-\eta\left([Z_{a,X}+b_{a,X}^r,Y]+[X,Z_{a,Y}+b_{a,Y}^r]\right)\\
=&\tg^*\left(\left(\nabla_{\bar X}^*\nabla^*_{\bar Y}\bar\eta-\nabla_{\bar Y}^*\nabla^*_{\bar X}\bar\eta\right)(\bar a)-[\bar X,\bar Y](\bar \eta(\bar a))\right)\\
&+Y(\eta(Z_{a,X}))-X(\eta(Z_{a,Y}))-\eta\left(Z_{b_{a,X},Y}+b_{b_{a,X},Y}^r-Z_{b_{a,Y},X}-b_{b_{a,Y},X}^r\right)\\
\overset{(3)}=&\tg^*\left(\bar\eta(-\nabla_{[\bar X,\bar Y]}\bar a
-\nabla_{\bar Y}\nabla_{\bar X}\bar a +\nabla_{\bar X}\nabla_{\bar Y}\bar a) \right)\overset{(3)}=0.
\end{align*}
Hence, $a^r (\eta([X,Y]))=0$ for all $a\in\Gamma(A)$ and since $G$ is
source-connected, this implies that
$\eta([X,Y])(g)=\eta([X,Y])(\s(g))$ for all $g\in G$. But since for $m\in M$, we
have
\[[X,Y](m)=[\bar X,\bar Y](m)\]
and $F_M$ is a subalgebroid of $TM$ by (1),
we find that $[X,Y](\s(g))\in F_G(\s(g))$ for all $g\in G$ and hence
$\eta([X,Y])(g)=\eta([X,Y])(\s(g))=0$.
Since $\eta$ was a $\tilde\tg$-section of $F_G^\circ$ and
$\tilde\tg$-sections 
of $F_G^\circ$ span $F_G^\circ$ on $G$, we have shown that 
$[X,Y]\in\Gamma(F_G)$ and the proof is complete.
\end{proof}

\begin{remark}
\begin{enumerate}
\item We have seen in this proof that Condition (2) implies the fact that $F_\core$ is a subalgebroid of $A$.
\item The same result has been shown independently in \cite{CrSaSt12}, using Lie groupoid and Lie algebroid cocycles,
in the special case where $F_M=TM$, i.e. where $F_G$ is a wide subgroupoid of $TG$.
\end{enumerate} 
\end{remark}

\subsection{Examples}
\begin{example}
Assume that $G$ is a Lie group (hence with $M=\{e\}$)
 with Lie algebra $\lie g$. Let $F_G$ be a multiplicative 
distribution. In this case, the core $F_\core=:\lie f$ 
is the fiber of $F_G$ over the identity and $F_{M}=0$. As a consequence, any partial $F_M$-connection 
on $\lie g/\lie f$ is trivial.
We check that all the conditions in Theorem \ref{the_connection} are automatically satisfied.

First of all, any element $\xi$ of $\lie g$ is 
$\nabla$-parallel. This implies that 
\[\left[\xi^r,\Gamma(F_G)\right]\subseteq \Gamma(F_G)\qquad \text{ for all }\qquad \xi\in\lie g,\]
i.e., $F_G$ is left-invariant, in agreement with \cite{Ortiz08,Jotz11a,Jotz11b}.

1), 3) and 4) are trivially satisfied and 2) is exactly the fact that $\lie f$ is an ideal 
in $\lie g$. This recovers the results proved in \cite{Ortiz08,Jotz11a,Jotz11b}.
\end{example}

\begin{example}\label{ex_group_action}
Let $G\rr M$ be a Lie groupoid with a smooth, free and proper
action of a Lie group $H$ by Lie groupoid automorphisms. 
Let $\mathcal V_G$ be the vertical space of the action, i.e.,
the smooth subbundle of $TG$ that is generated by the infinitesimal 
vector fields $\xi_G$, for all $\xi\in\lie h$, where $\lie h$ is the Lie algebra
of $H$.  The involutive subbundle $\mathcal V_G$ is easily seen to be 
multiplicative (see for instance \cite{Jotz11b}).

The action restricts to a free and proper action
of $H$ on $M$, and it is easy to check that
$\mathcal V_G\cap TM=\mathcal V_M$
is the vertical vector space of the action of $H$ on $M$.
Furthermore, $\mathcal V_G\cap T^\s G=\mathcal V_G\cap T^\tg G=0^{TG}$
and we get $\mathcal V^\core=0^A$.

The infinitesimal vector fields $(\xi_G,\xi_M)$, $\xi\in\lie h$, are multiplicative (in the sense of 
\cite{MaXu98} for instance). We get hence from \cite{MaXu98}
that the Lie bracket 
$[a^r, \xi_G]$ is right-invariant 
for any $\xi\in\lie h$ and $a\in\Gamma(A)$.
We obtain  a map (see also \cite{MaXu98})
\begin{align*}
\begin{array}{ccl}
\lie h\times\Gamma(A)&\to& \Gamma(A)\\
(\xi,a)&\mapsto&[\xi_G,a^r]\an{M}
\end{array},
\end{align*}
and  we recover the connection
\begin{align*}
\begin{array}{cccl}
\nabla:&\Gamma(\mathcal V_M)\times\Gamma(A)&\to& \Gamma(A)
\end{array}
\end{align*}
defined by $\nabla_{\xi_M}a=[\xi_G,a^r]\an{M}$
for all $\xi\in\lie h$ and $a\in\Gamma(A)$. This connection is obviously flat and satisfies 
all the conditions in Theorem \ref{the_connection}.
\end{example}

\begin{example}
Let $(G\rr M, J_G)$ be a complex Lie groupoid, i.e., a
Lie groupoid endowed with a complex structure
$J_G$
that is multiplicative in the sense that the map
\begin{equation*}
\begin{xy}
\xymatrix{
TG\ar[r]^{J_G}\ar@<.6ex>^{T\s}[d]\ar@<-.6ex>_{T\tg}[d]&TG\ar@<.6ex>^{T\s}[d]\ar@<-.6ex>_{T\tg}[d]\\
TM\ar[r]_{J_M}&TM
}
\end{xy}
\end{equation*}
is a Lie groupoid morphism over some map $J_M$.
Since ${J_G}^2=-\Id_{TG}$, we conclude that ${J_M}^2=-\Id_{TM}$
and the Nijenhuis condition for $J_M$
is easy to prove using $\s$-related vector fields.
The map $J_G$ restricts also to a map $j_A$ on the core $A$, i.e., a
fiberwise complex structure
that satisfies also a Nijenhuis condition. (This can be seen by noting
that the Nijenhuis tensor of $J_G$
restricts to right-invariant vector fields.)

The
subbundles
$T^{1,0}G=E_i$ and $T^{0,1}G=E_{-i}$ of $TG\otimes\mathbb C$
are multiplicative and involutive
with bases
$T^{1,0}M$ and $T^{0,1}M$  and cores
$A^{1,0}$ and $A^{0,1}$.
The quotient $(A\otimes\mathbb C)/A^{1,0}$
is isomorphic as a vector bundle to $A^{0,1}$
and a
 straightforward computation shows that the connection
that we get from the multiplicative complex foliation $T^{1,0}G$
is exactly
the connection
$\nabla:\Gamma(T^{1,0}M)\times \Gamma(A^{0,1})\to \Gamma(A^{0,1})$
as in Lemma 4.7 of \cite{LaStXu08}.

Since the parallel sections of the connection in this Lemma are exactly
the holomorphic sections of $A^{0,1}$,
one can reconstruct the map
$J_A:TA\to TA$ defined by $J_A=\sigma\inv\circ A(J_G)\circ \sigma$ \footnote{
To avoid confusions, we write in this example $\sigma:TA\to A(TG)$
for the canonical flip map.} as
in \cite{LaStXu09}
by requiring that $J_A(Ta)=Ta\circ J_M$
for all parallel sections $a\in\Gamma(A)$, and $J_A(\hat
b)=\widehat{j_A(b)}\circ J_M$
for all sections $b\in\Gamma(A)$.

By  Lemma 4.7 in \cite{LaStXu08} and the integration results in
\cite{LaStXu09}, the complex structure
$J_G$ is hence equivalent to the data $J_M, j_A$ and this connection with its properties.
This is
in agreement with the results that we will prove in section
\ref{connection-as-infinitesimal-data}.
\end{example}

\vspace*{0.5cm}

\section{Foliated algebroids}\label{fol_alg}

In this section we study Lie algebroids equipped with foliations compatible with both 
Lie algebroid structures $TA\to TM$ and $TA\to A$ on $TA$. This is the first step towards an infinitesimal 
description of multiplicative foliations.

\subsection{Definition and properties}
\begin{definition}
Let $A\to M$ be a Lie algebroid. A subbundle 
$F_A\subseteq TA$ is called \textbf{morphic} if it is a Lie subalgebroid 
of $TA\to TM$ over some subbundle $F_M\subseteq TM$.

If $F_A$ is involutive and morphic, then the pair $(A, F_A)$ is referred to as a \textbf{foliated Lie algebroid}.
\end{definition}

Consider a foliated Lie algebroid $(A,F_A)$.
We have the $\mathcal{VB}$-algebroid \begin{align*}
\begin{xy}
\xymatrix{
F_A\ar[r]^{Tq_A\an{F_A}}\ar[d]_{p_A\an{F_A}}&F_M\ar[d]^{p_M\an{F_M}}\\
A\ar[r]_{q_A}&M
}
\end{xy}
\end{align*}
and the Lie algebroid morphisms
\begin{align*}
\begin{xy}
\xymatrix{
F_A\ar[d]_{Tq_A\an{F_A}}\ar[r]^{\iota_{F_A}}&TA\ar[d]^{Tq_A}\\
F_M\ar[r]_{\iota_{F_M}}&TM
}
\end{xy}\qquad 
\begin{xy}
\xymatrix{
F_A\ar[d]_{p_A\an{F_A}}\ar[r]^{\iota_{F_A}}&TA\ar[d]^{p_A}\\
A\ar[r]_{\id_A}&A
}
\end{xy}.
\end{align*}
We have 
\begin{equation}\label{splitting_F_A}
F_A\left(0^A(m)\right)=T_m0^A(F_M(m))\oplus\ker(Tq_A\an{F_A(0^A(m))})
\end{equation}
for all $m\in M$.
 If 
$v_{0^A(m)}\in F_A(0^A_m)$ is such that 
$Tq_A(v_{0^A(m)})=0_m\in F_M(m)$, then 
\[v_{0^A(m)}=\left.\frac{d}{dt}\right\an{t=0}
0_m^A+ta_m\]
for some $a_m\in A_m$.
Set \[F_\core(m)
:=\left\{a_m\in A_m\left|\left.\frac{d}{dt}\right\an{t=0}
0_m^A+ta_m \in F_A\left(0^A_m\right)\right.\right\}\]
for all $m\in M$. Equation \eqref{splitting_F_A} shows that this defines
a vector subbundle of $A$. 

Furthermore, using the fact that $F_A$ is closed under 
the addition in $TA\to TM$, one can check that
\begin{equation}\label{kernel_core}
\ker\left(Tq_A\an{F_A(b_m)}\right)
=\left\{\left.\left.\frac{d}{dt}\right\an{t=0}
b_m+ta_m \right| a_m\in F_\core(m)\right\}
\end{equation}
for all $b_m\in A_m$.
The vector bundle $F_\core\simeq \ker(Tq_A)\cap \ker(p_A)$ is the \textbf{core} of $F_A$.

Recall that the core sections of $TA\to A$ are the sections $a^\uparrow$ as in \eqref{core_uparrow} for $a\in \Gamma(A)$, 
and the core sections of $TA\to TM$ are the sections $\hat a$ as in \eqref{eq:coreTqA}. The core sections of $TA\to A$
that have image in $F_A$ are hence exactly the core sections $a^\uparrow$ for $a\in\Gamma(F_\core)$.
Now if 
\[\hat a(v_m)=T_m0^Av_m+_{p_A}\left.\frac{d}{dt}\right\an{t=0}ta(m)\in F_A\]
 then we have $v_m=Tq_A(\hat a(v_m))\in F_M(m)$
and since $F_A\to F_M$ is a vector subbundle of $TA\to TM$,
the zero element $T_m0^Av_m$ belongs to the fiber of $F_A$ over $v_m$.
This leads to $\left.\frac{d}{dt}\right\an{t=0}ta(m)\in F_A$, using the addition
in the fiber of $F_A$ over $v_m$, and hence
$a(m)\in F_\core(m)$.
\emph{This shows that the core sections $\hat a$ of $TA\to TM$ that restrict to core sections
of $F_A\to F_M$ are the restrictions 
$\hat a\an{F_M}$ for $a\in\Gamma(F_\core)$.}

The vector bundles  $F_A\to A$ and $F_A\to F_M$ are spanned by their linear and core sections.
\begin{lemma}\label{lem_a_uparrow_preserves_FA}
Let $(A,F_A)$ be  a foliated Lie algebroid.
A section $a\in\Gamma(A)$ is such that $Ta\an{F_M}\in\Gamma_{F_M}(F_A)$ if and only if $\left[a^\uparrow,\Gamma(F_A)\right]
\subseteq \Gamma(F_A)$.
\end{lemma}

\begin{proof} 
Choose any linear section $X$ of $F_A\to A$
over a non-vanishing section $\bar X\in\Gamma(F_M)$. Since 
$\bar X(m)\neq 0$ for all $m$, we get that for each $c_m\in A_m$, the vector 
$X(c_m)$ can be written $X(c_m)=T_{m}c(\bar X(m))$ for some section 
$c\in \Gamma(A)$ with $c(m)=c_m$.
We write for simplicity $\phi_{\cdot}^a$ for the flow 
of $a^\uparrow$, i.e., $\phi_t^a(c')=c'+ta(q_A(c'))$ for all $c'\in A$ and $t\in\R$, and 
$\bar\phi_{\cdot}$ for the flow of $\bar X$.
We have,  writing $v_m=\bar X(m)$,
\begin{align}
T_{c_m}\phi^a_tX(c_m)&=T_{m}(\phi^a_t\circ c)(\bar X(m))\nonumber\\
&=\left.\frac{d}{ds}\right\an{s=0}c\left(\bar\phi_s(m)\right)+ta\left(\bar\phi_s(m)\right)\nonumber\\
&=T_{m}c(v_m)+tT_ma(v_m)\nonumber\\
&=X(c_m)+tT_ma(v_m).\label{eq_Ta}
\end{align}
If  $\left[a^\uparrow, \Gamma(F_A)\right]\subseteq \Gamma(F_A)$,
the flow $\phi^a$ of $a^\uparrow$ preserves $F_A$.
Hence,  $T_{c_m}\phi^a_tX(c_m)\in F_A(c_m+ta(m))$
is a vector of $F_A\to F_M$ over $\bar X(m)$ and since 
$X(c_m)$ is also a vector of $F_A\to F_M$ over $\bar X(m)$, we find that 
$tT_ma(v_m)\in F_A(ta(m))$ over $\bar X(m)$ for all $t\in\R$. 
In particular we  get  $Ta\an{F_M}\in\Gamma_{F_M}(F_A)$.

Conversely, if $Ta\an{F_M}\in\Gamma_{F_M}(F_A)$, then 
$tT_ma(v_m)\in F_A(ta(m))$ for all $t\in\R$ and $v_m\in F_M$. Since $X(c_m)$
is an element of $F_A$ for all $c_m\in A_m$, we find hence $T_{c_m}\phi^a_tX(c_m)\in F_A(c_m+ta(m))$
and  $[a^\uparrow, X]\in\Gamma(F_A)$. This implies 
the inclusion 
$\left[a^\uparrow, \Gamma(F_A)\right]\subseteq \Gamma(F_A)$
since the bracket of two core sections is $0$ and $F_A\to A$ is spanned by its linear and its core sections.
\end{proof}

\begin{proposition}\label{core_and_base_alg}
Let $(A,F_A)$ be  a foliated Lie algebroid.
Then  the core $F_\core\to M$ of $F_A$ is a subalgebroid of $A$ and
the base $F_M\to M$ is an involutive subbundle of $TM$. 
\end{proposition}
\begin{proof}
Choose $a,b\in\Gamma(F_\core)$. Since $F_A\to A$ is involutive and $a^\uparrow\in\Gamma(F_A)$, we know 
that $\left[a^\uparrow,\Gamma(F_A)\right]\subseteq \Gamma(F_A)$.
By the preceding lemma, this implies that $Ta\an{F_M}$ is a section of $F_A\to F_M$.
Since $\hat b\an{F_M}$ is also a section of the subalgebroid $F_A\to F_M$ of $TA\to TM$, we find 
that $\widehat{[a,b]}\an{F_M}=\left[Ta\an{F_M},\hat b\an{F_M}\right]
\in\Gamma_{F_M}(F_A)$ and hence 
$[a,b]\in\Gamma(F_\core)$.

The second statement is immediate by Lemma \ref{linear_core_of_TA_A} and the fact that
$q_A$ is a surjective submersion.
\end{proof}

The following lemma will  also be useful later.
\begin{lemma}\label{linear_over_0}
Let $(A, F_A)$ be a foliated Lie algebroid. Let $\left(X,0^{TM}\right)$ be a linear section of $TA$ covering 
the zero section.
Then we have 
\[X\in\Gamma(F_A)
\Leftrightarrow
D_Xa\in\Gamma(F_\core) \text{ for any } a\in\Gamma(A).\]
\end{lemma}

\begin{proof} Choose $a_m\in A$.
Since $T_{a_m}q_A(X(a_m))=0^{TM}_m$, we can write
$X(a_m)=(c(a_m))^\uparrow(a_m)$ for some 
$c(a_m)\in\Gamma(A)$.
We have then, for any section $\varphi\in\Gamma((F_\core)^\circ)$, 
where $(F_\core)^\circ$ is the annihilator of $F_\core$ in $A^*$,
and any section $a\in\Gamma(A)$ such that $a(m)=a_m$:
\begin{align*}
\langle\varphi(m),D_Xa(m)\rangle&=\left.\frac{d}{dt}\right\an{t=0}\langle \varphi(m),\phi^X_{-t}(a_m)\rangle
=-X(a_m)(l_\varphi)\\
&=\left.\frac{d}{dt}\right\an{t=0}\langle \varphi(m),a_m-tc(a_m)(m)\rangle
=-\langle \varphi(m),c(a_m)(m)\rangle.
\end{align*}
Hence,  $(D_Xa)(m)\in F_\core(m)$  if and only if $c(a_m)(m)\in F_\core(m)$.
If $D_Xa\in\Gamma(F_\core)$ for all $a\in\Gamma(A)$, we find hence $X\in\Gamma(F_A)$ 
by \eqref{kernel_core}, and conversely, if 
$X\in\Gamma(F_A)$, we get $D_Xa\in\Gamma(F_\core)$ for all $a\in\Gamma(A)$.
\end{proof}

\subsection{The Lie algebroid of a multiplicative foliation}

The following construction can be found in \cite{Ortiz08t} in the more general setting of multiplicative Dirac structures.

Let $F_G$ be a multiplicative subbundle of $TG$ with space of units $F_M\subseteq TM$. 
Since $F_G\subseteq TG$ is a Lie subgroupoid, we can apply the Lie functor, 
leading to a Lie subalgebroid $A(F_G)\subseteq A(TG)$
over $F_M\subseteq TM$.

 As we have seen in Subsection \ref{tangentalgebroids}, the canonical involution $J_G:TTG\rmap TTG$ 
restricts to an isomorphism of double vector bundles $j_G:TA\rmap A(TG)$ inducing the identity map 
on both the side bundles and the core. Since $j_G:TA\rmap A(TG)$ is an isomorphism of Lie algebroids 
over $TM$, we conclude that 
 \[F_A:= j^{-1}_G(A(F_G))\subseteq TA\]
is a Lie algebroid over $F_M\subseteq TM$. 
Since 
\begin{align*}
\begin{xy}
\xymatrix{
F_G\ar[r]^{p_G}\ar@<.6ex>^{T\s}[d]\ar@<-.6ex>_{T\tg}[d]& G\ar@<.6ex>^{\s}[d]\ar@<-.6ex>_{\tg}[d] \\
F_M\ar[r]_{p_M}&M}
\end{xy}
\end{align*}
is a $\mathcal{VB}$-subgroupoid of 
\begin{align*}
\begin{xy}
\xymatrix{
TG\ar[r]^{p_G}\ar@<.6ex>^{T\s}[d]\ar@<-.6ex>_{T\tg}[d]& G\ar@<.6ex>^{\s}[d]\ar@<-.6ex>_{\tg}[d] \\
TM\ar[r]_{p_M}&M}
\end{xy},
\end{align*}
the Lie algebroid
\begin{align*}
\begin{xy}
\xymatrix{
A(F_G)\ar[r]^{A(p_G)}\ar[d]&A\ar[d] \\
F_M\ar[r]_{p_M}&M}
\end{xy}
\end{align*}
is a $\mathcal{VB}$-subalgebroid of 
\begin{align*}
\begin{xy}
\xymatrix{
A(TG)\ar[r]^{A(p_G)}\ar[d]& A\ar[d] \\
TM\ar[r]_{p_M}&M}
\end{xy}
\end{align*}
(\cite{BuCadelHo11}), and $F_A\to A$ is also a subbundle of $TA\to A$. 

Before we study the sections of $A(F_G)\to F_M$, we show that the vector subbundle 
$F_A\subseteq TA$ is involutive. The space of sections of $F_A\to A$ is spanned by the 
core sections $a^\uparrow$ for $a\in\Gamma(F_\core)$ and 
linear sections $(\tilde X,\bar X)=(j_G\inv(A(X)),\bar X)$ coming from star 
sections $X\overset{\star}\sim_\s\bar X$ of $F_G\to F_M$. Hence, the involutivity of $F_A$ follows immediately
from the involutivity of $F_G\subseteq TG$ and the considerations following Lemma \ref{linear_core_of_TA_A}.

\medskip

The Lie algebroid $A(F_G)\to F_M$
is a subalgebroid of $A(TG)\to TM$. We want to find the  linear sections 
of $A(TG)$ that restrict to linear sections of $A(F_G)$.
If $a\in\Gamma(F_\core)$, then 
$\tilde a^r$ is easily seen to be tangent to $F_G$ on $F_M$.
Choose an arbitrary section $a\in\Gamma(A)$. We want to find a condition for $\beta_a\an{F_M}$\footnote{Recall 
the definitions of $\tilde a$ and $\beta_a$ from \eqref{a_tilde} and \eqref{beta_a}.}
to be a section of $A(F_G)$.
We have to find 
\[\beta_a(v)\in T_{v}F_G\]
for all $v\in F_M$. A sufficient condition is $[a^r,\Gamma(F_G)]\subseteq \Gamma(F_G)$.
We will show in the remainder of this subsection that this condition is also necessary.

\bigskip

Now we define a more natural, partial 
$F_M$-connection $\nabla^\s$ on $A/F_\core$, that 
 will help us to understand which linear sections of $A(TG)\to TM$
restrict to linear sections of $A(F_G)\to F_M$.
Choose $a\in\Gamma(A)$, $\bar X\in\Gamma(F_M)$ and set
\begin{equation}\label{nabla_s}
\nabla^\s_{\bar X}\bar a:=\overline{[X,a^r]\an{M}}
\end{equation}
for any section $X\in\Gamma(F_G)$ such that $X\overset{\star}\sim_\s\bar X$.
We show the following proposition.
\begin{proposition}
Let $(G\rr M,F_G)$ be a foliated Lie groupoid. Then 
$\nabla^\s$  as in \eqref{nabla_s} is a well-defined flat partial $F_M$-connection on $A/F_\core$.
\end{proposition}

\begin{proof}
For any section $a\in\Gamma(A)$, we have $a^r\sim_\s0$. Since 
$X\overset{\star}\sim_\s\bar X$ implies in particular
$X\sim_\s\bar X$, we find 
$[X,a^r]\sim_\s 0$ and $[X,a^r]\an{M}$ is thus a section of $A$.

We show that $\overline{[X,a^r]\an{M}}\in\Gamma(A/F_\core)$ does not depend 
on the choice of $X\in\Gamma(F_G)$. 
If $X'\in\Gamma(F_G)$, $X'\overset{\star}\sim_\s\bar X$, then
we have $X-X'\in\Gamma(F_G)$ and 
$X-X'\overset{\star}\sim_\s0$. Hence, if $a_1,\ldots, a_r$
is a local frame for $F_\core$ on an open set $U\subseteq M$, 
we can write 
$X-X'=\sum_{i=1}^rf_ia_i^r$
on $\tilde U:=\tg\inv(U)$, with $f_1,\ldots,f_r\in C^\infty(\tilde U)$
such that $f_i\an{U}=0$ for $i=1,\ldots,r$.
We have then 
\[[X-X',a^r]=-\sum_{i=1}^ra^r(f_i)a_i^r+\sum_{i=1}^rf_i[a_i,a]^r.
\]
But since the second term of this sum vanishes on $U$, we conclude that
$[X-X',a^r]\an{M}\in\Gamma(F_\core)$. This proves that $\nabla^\s$ is well-defined.

Assume now that $X\overset{\star}\sim_\s \bar{X}$ and choose $f\in C^\infty(M)$. Then we have 
$(\s^*f)\cdot X\overset{\star}\sim_\s f\cdot \bar X$ and
\begin{align*}
[\s^*f\cdot X, a^r]=\s^*f\cdot[ X, a^r]-a^r(\s^* f)X=\s^*f\cdot[ X, a^r],
\end{align*}
which shows
\[\nabla^\s_{f\cdot \bar X}\bar a=f\cdot\nabla^\s_{\bar X}\bar a.
\]
Analogously, $(f\cdot a)^r=\tg^*f\cdot a^r$
and, since $X\an{M}=\bar X\in\mx(M)$,
\begin{align*}
[X, \tg^*f\cdot a^r]=\tg^*f\cdot[ X, a^r]+X(\tg^*f)\cdot a^r
\end{align*}
restricts to 
\[f\cdot[ X, a^r]\an{M}+\bar X(f)\cdot a
\]
on $M$.

We show finally that $\nabla^\s$ is flat.
Choose $X\overset{\star}\sim_\s\bar X$ and $Y\overset{\star}\sim_\s\bar Y$.
Then, since 
\[[X,a^r]-\left([X,a^r]\an{M}\right)^r\]
vanishes on $M$
and $Y$ is  tangent to $M$ on $M$, we find that
\[\bigl[Y, [X,a^r]-\left([X,a^r]\an{M}\right)^r\bigr]\]
vanishes on $M$.
Thus, we have 
\[\bigl[Y, [X,a^r]\bigr]\an{M}=\bigl[Y,\left([X,a^r]\an{M}\right)^r\bigr]\an{M}\]
and the flatness of $\nabla^\s$ follows easily.
\end{proof}
Recall that we say that $a\in\Gamma(A)$ is $\nabla$- (or $\nabla^\s$-)parallel if $\bar a$ is $\nabla$-parallel, 
since this doesn't depend on the chosen representative.
\begin{theorem}\label{char_of_nabla_s_parallel_sections}
Let $(G\rr M,F_G)$ be a foliated Lie groupoid and choose $a\in\Gamma(A)$. Then the following are equivalent:
\begin{enumerate}
\item \[[a^r,\Gamma(F_G)]\subseteq\Gamma(F_G),\]
\item $a$ is $\nabla$-parallel, where $\nabla$ is the connection in Theorem \ref{the_connection},
\item $a$ is $\nabla^\s$-parallel,
\item $\beta_a\an{F_M}$ is a section of $A(F_G)$,
\item $[X,a^r]\an{M}\in\Gamma(F_\core)$
for any section $X\in\Gamma(F_G)$ such that $X\overset{\star}\sim_\s\bar X$.
\end{enumerate}
\end{theorem}

This theorem shows in particular that the $\nabla$- (or equivalently $\nabla^\s$-) parallel
sections of $A$ are \emph{exactly} the sections of $A$ such that
$\beta_a\an{F_M}$ is a section of $A(F_G)$. 

Since the two connections are flat and their parallel sections coincide, we get the following corollary.

\begin{corollary}\label{char_of_nabla_s_parallel_sections_cor}
Let $(G\rr M,F_G)$ be a foliated Lie groupoid.
Then $\nabla=\nabla^\s$.
\end{corollary}

\begin{remark}
Note that in Example \ref{ex_group_action}, the two connections can immediatly be seen
to be the same by definition, since the multiplicative vector fields $\xi_G$, $\xi\in\lie h$
are $\tg$-sections and star sections at the same time.
\end{remark}

\begin{proof}[Proof of Theorem \ref{char_of_nabla_s_parallel_sections}]
We already know by Theorem \ref{the_connection} that $(1)$ and $(2)$ are equivalent.
By definition of $\nabla^\s$, $(3)$ and $(5)$ are equivalent.

We show that $a$ is $\nabla$-parallel if and only if  $a$ is $\nabla^\s$-parallel.
If $a$ is $\nabla$-parallel, then 
$[a^r,\Gamma(F_G)]\subseteq \Gamma(F_G)$. This yields immediately
$[X,a^r]\an{M}\in\Gamma(F_\core)$ for any star-section $X\overset{\star}\sim_\s \bar X$
of $F_G$. Conversely, assume that $a$ is $\nabla^\s$-parallel. Since $\nabla$ is flat, $a$
can be written locally as  a sum $a=\sum_{i=1}^nf_ia_i$
with $f_i\in C^\infty(M)$, $i=1,\ldots,n$ and $\nabla$-parallel sections 
$a_1,\ldots,a_n$.
We have then for any section $\bar X$ of $F_M$: 
\[\nabla_{\bar X}\bar a=\sum_{i=1}^n \bar X(f_i)\bar{a_i}.
\]
But since the sections $a_i$, $i=1,\ldots,n$ are then also $\nabla^\s$-parallel by the considerations above, we have also
\[\sum_{i=1}^n \bar X(f_i)\bar{a_i}=\nabla^\s_{\bar X}\bar a=0,
\]
which shows that $a$ is $\nabla$-parallel.

We show that $(1)$ implies $(4)$. If $\left[a^r,\Gamma(F_G)\right]\subseteq \Gamma(F_G)$, then 
the flow $L_{\Exp(\cdot a)}$ of $a^r$ leaves $F_G$ invariant by Corollary \ref{appendix_cor}.
Hence, we have $T_m\Exp(ta)v_m\in F_G(\Exp(ta)(m))$
for all $v_m\in F_M(m)$. This yields 
$\beta_a(v_m)=\left.\frac{d}{dt}\right\an{t=0}T_m\Exp(ta)v_m\in T_{v_m}^{T\s}F_G=A_{v_m}(F_G)$.

We show that $(4)$ implies $(5)$.
Recall that the flow of ${\beta_a}^r$ is equal
to $TL_{\Exp(\cdot a)}$. Hence, if $\beta_a\an{F_M}\in A(F_G)$, we have 
that \[T_mL_{\Exp(sa)}(v_m)\in F_G(\Exp(sa)(m))\]
for all $|s|\in\R$ small enough and $v_m\in F_M$.
We have also 
\[T_mL_{\Exp(sa)}(u_m)=0_{\Exp(sa)(m)}\star u_m\in F_G(\Exp(sa)(m))\]
for all $u_m\in F^\tg$. This yields
\[T_mL_{\Exp(sa)}F_G(m)\subseteq F_G(\Exp(sa)(m))\]
for all $m\in M$. Since both sides are vector spaces of the same dimension,
we get an equality
and in particular
\[F_G(m)=T_{\Exp(sa)(m)}L_{\Exp(-sa)}F_G(\Exp(sa)(m)).\]
 Choose now any star section $X\overset{\star}\sim_\s\bar X$.
We have for any $m\in M$:
\[[a^r,X](m)=\left.\frac{d}{ds}\right\an{s=0}
T_{\Exp(sa)(m)}L_{\Exp(-sa)}X(\Exp(sa)(m)).\]
Since $T_{\Exp(sa)(m)}L_{\Exp(-sa)}X(\Exp(sa)(m))\in F_G(m)$
for all $s$ small enough, we find that $[a^r,X](m)\in F_\core(m)$.
\end{proof}

\subsection{Integration of foliated algebroids}
The main theorem of this subsection is the following. 
\begin{theorem}\label{integration_of_fol_alg}
Let $(G\rr M, F_G)$ be a foliated groupoid. Then 
$(A, F_A=j_G\inv(A(F_G)))$ is a foliated algebroid.

Conversely, let $(A,F_A)$ be a foliated Lie algebroid. Assume that
$A$ integrates to a source simply connected Lie groupoid $G\rr M$. Then there is a unique multiplicative 
foliation $F_G$ on $G$ such that $F_A=j_G\inv(A(F_G))$.
\end{theorem}

We will use a result of \cite{BuCadelHo11}, which 
states that a $\mathcal{VB}$-algebroid 
\begin{align*}
\begin{xy}
\xymatrix{
 E\ar[r]^{q_{ E}^h}\ar[d]_{q_{ E}^v}&A\ar[d]^{q_A}\\
B\ar[r]_{q_B}&M
}
\end{xy}
\end{align*}
integrates to a $\mathcal{VB}$-groupoid
\begin{align*}\begin{xy}
\xymatrix{
G(E)\ar[r]^{q_{G(E)}}\ar@<.6ex>[d]\ar@<-.6ex>[d]& G(A)\ar@<.6ex>[d]\ar@<-.6ex>[d] \\
B\ar[r]_{q_B}&M}
\end{xy}.
\end{align*}
Furthermore, if 
$E'\hookrightarrow E$, $B'\hookrightarrow B$ is a $\mathcal{VB}$-subalgebroid 
with the same horizontal base $A$,
\begin{align*}
\begin{xy}
\xymatrix{
E'\ar[r]^{q_{ E}^h}\ar[d]_{q_{ E}^v}&A\ar[d]^{q_A}\\
B'\ar[r]_{q_{B'}}&M
}
\end{xy}, 
\end{align*}
then $E'\to B'$ integrates to an embedded $\mathcal{VB}$-subgroupoid 
$G(E')\hookrightarrow G(E)$ over $B'\hookrightarrow B$,
 \begin{align*}\begin{xy}
\xymatrix{
G(E')\ar[r]^{q_{G(E)}}\ar@<.6ex>[d]\ar@<-.6ex>[d]& G(A)\ar@<.6ex>[d]\ar@<-.6ex>[d] \\
B'\ar[r]_{q_M}&M}
\end{xy}.
\end{align*}
This is done in \cite{BuCadelHo11}
using the characterization of vector bundles via homogeneous structures  (see \cite{GrRo09}).

\begin{proof}[Proof of Theorem \ref{integration_of_fol_alg}]
The first part has been shown in the previous subsection.

\medskip

Let $(A,F_A)$ be a foliated Lie algebroid. Then the core $F_\core\subseteq A$ and the base
$F_M\subseteq TM$ are subalgebroids by Proposition \ref{core_and_base_alg}.
The $\mathcal{VB}$-subgroupoid of $(TG,G;TM,M)$
integrating the subalgebroid $j_G(F_A)\to F_M$
of $A(TG)\to TM$  is a multiplicative subbundle $F_G\rr F_M$
of $TG\rr TM$ with core $F_\core$. By Theorem \ref{inv_crit}, $F_G$ is involutive.
\end{proof}

\subsection{The connection associated to a foliated algebroid}
Let $A$ be a Lie algebroid endowed with an involutive 
subbundle $F_A\subseteq TA$. 
We show that if $F_A\to F_M$ is also morphic, that is, a subalgebroid
of the tangent algebroid $TA\to TM$,  
then the Bott connection 
\[\nabla^{F_A}:\Gamma(F_A)\times \Gamma(A/F_A)\to\Gamma(A/F_A)
\]
induces
a natural partial $F_M$-connection 
on $A/F_\core$ with the same properties as the connection in Theorem \ref{the_connection}.

For $\bar X$ in $\Gamma(F_M)$, we choose any linear section $X\in\Gamma_A(F_A)$
over $\bar X$. Then, for any
core section $a^\uparrow$ of $TA\to A$, 
we have 
\[ \left[X,a^\uparrow\right]=D_Xa^\uparrow\]
for some section $D_Xa$ of $A$, see Lemma \ref{linear_core_of_TA_A}.
Set \begin{equation}\label{nabla_A}
\nabla^{A}_{\bar X}\bar a=\overline{D_Xa} \in \Gamma(A/F_\core).
\end{equation}

We show the following proposition.
\begin{proposition}\label{nabla_A_welldef}
Let $(A,F_A)$ be a foliated Lie algebroid. The map 
\[\nabla^{A}:\Gamma(F_M)\times\Gamma(A/F_\core)\to \Gamma(A/F_\core)
\] as defined in \eqref{nabla_A}
is a well-defined flat partial $F_M$-connection on $A/F_\core$.
\end{proposition}

\begin{proof}
Choose $\bar X\in\Gamma(F_M)$
and two linear sections $X,Y\in\Gamma_A(F_A)$
covering $\bar X$.
Then $(X-Y,0)$ is a vector bundle 
morphism and $Z:=X-Y\in\Gamma_A(F_A)$. 
By Lemma \ref{linear_over_0}, 
we get $D_Za\in\Gamma(F_\core)$ for any $a\in\Gamma(A)$
and since 
\[(D_Za)^\uparrow=[X-Y,a^\uparrow]=(D_Xa-D_Ya)^\uparrow,
\] this shows that $\nabla^A$ is well-defined. The properties 
of a connection can be checked using the equality
$\left(f\cdot a\right)^\uparrow=q_A^*f\cdot a^\uparrow$
for $a\in\Gamma(A)$ and $f\in C^\infty(M)$. 

The flatness of $\nabla^A$ follows immediately from the flatness of the Bott connection $\nabla^{F_A}$.
\end{proof}

\begin{example} 
Let $H$ be a connected Lie group with Lie algebra $\mathfrak{h}$. Assume that $H$ acts on a Lie algebroid $A\rmap M$ in 
a free and proper manner, by Lie algebroid automorphisms. That is, for
all $h\in H$, the diffeomorphism 
$\Phi_h$ is a Lie algebroid morphism 
over $\phi_h:M\to M$. Consider 
the vertical spaces $\mathcal{V}_A$, $\mathcal{V}_M$ defined as follows
\[\mathcal{V}_A(a)=\{\xi_A(a)\mid \xi \in\mathfrak{h}\}, \qquad \mathcal{V}_M(m)=\{\xi_M(m)\mid \xi \in\mathfrak{h}\},\]
for  $a\in A$ and $m\in M$. 
We check that $\mathcal{V}_A$ inherits a Lie algebroid structure over $\mathcal{V}_M$ making the pair 
$(A,\mathcal{V}_A)$ into a foliated algebroid with core zero. Choose
$\xi_A(a_m)\in\V_A(a_m)$ for some $\xi\in\lie g$, then
$T_{a_m}q_A\xi_A(a_m)=\xi_M(m)\in\V_M(m)$. Hence, if
$T_{a_m}q_A\xi_A(a_m)=0$, then $\xi=0$ since the action is supposed to
be free and we find $F_\core=0$. 
Notice that, since the action is by algebroid 
automorphisms,  the infinitesimal generators $\xi_A$ are in fact morphic vector fields covering $\xi_M$.

We have, writing $\rho(a_m)=\dot c(0)$,
\begin{align*}
\rho_{TA}(\xi_A(a_m))
&=J_M\left(T\rho(\xi_A(a_m)
\right)=J_M\left(\left.\frac{d}{dt}\right\an{t=0} (\rho\circ\Phi_{\exp(t\xi)})(a_m)
\right)\\
&=J_M\left(\left.\frac{d}{dt}\right\an{t=0} (T\phi_{\exp(t\xi)}\circ \rho)(a_m)
\right)=\left.\frac{d}{ds}\right\an{s=0}\xi_M(c(s))\in T_{\xi_M(m)}\V_M.
\end{align*}
If $a\in\Gamma(A)$ is such that $T_ma(\xi_M(m))\in\V_A(a_m)$ for some
$\xi\in\lie g$ and $m\in M$, then there exists 
$\eta\in\lie g$ such that $T_ma(\xi_M(m))=\eta_A(a(m))$. But applying
$T_{a_m}q_A$ to both sides of this equality
yields then $\xi_M(m)=\eta_M(m)$, which leads to $\xi=\eta$, since the
action is free, and hence $T_ma(\xi_M(m))=\xi_A(a(m))$.

Here, 
the induced partial $\mathcal{V}_M$-connection on $A$ is given by $\nabla^A_{\xi_M}a\in\Gamma(A)$ where
\[[\xi_A,a^ \uparrow]=(\nabla^A_{\xi_M}a)^\uparrow\]
for any $a\in\Gamma(A)$.
If $a\in\Gamma(A)$ is $\nabla$-parallel, then 
we find $[\xi_A,a^ \uparrow]=0$ for all $\xi\in\lie g$ and 
hence 
\begin{align*}
0=&\left.\frac{d}{dt}\right\an{t=0}\left.\frac{d}{ds}\right\an{s=0}\Phi_{\exp(s\xi)}(b_m+ta(m))-ta(\phi_{\exp(s\xi)})\\
=& \left.\frac{d}{dt}\right\an{t=0}\xi_A(b_m)
+t\cdot \left.\frac{d}{ds}\right\an{s=0}\left(\Phi_{\exp(s\xi)}(a(m))-a(\phi_{\exp(s\xi)})\right)\\
\end{align*}
for all $b_m\in A$.
This leads to $\xi_A(a(m))=T_ma\xi_M(m)$ and 
hence
\begin{equation}\label{par_sec_V_A}
Ta(\V_M)\subseteq \V_A.
\end{equation}
The parallel sections of $\nabla^A$ are hence exactly the sections of $A$ satisfying \eqref{par_sec_V_A}. 
This will be shown in the general situation in the next lemma.

If $a\in\Gamma(A)$ is $\nabla^A$-parallel, we find easily that $a(\phi_h(m))=\Phi_h(a(m))$ for all 
$h\in H$ (recall that $H$ is assumed to be connected) and $m\in M$.
Hence, if $a,b\in\Gamma(A)$ are such that $Ta(\V_M)\subseteq \V_A$ and 
 $Tb(\V_M)\subseteq \V_A$, we have 
$a\circ \phi_h=\Phi_h\circ a$ and $b\circ\phi_h=\Phi_h\circ b$ for all $h\in H$.
Since the action is by Lie algebroid morphisms, we get then
$[a,b]\circ \phi_h=\Phi_h\circ [a,b]$, which implies 
$T[a,b](\V_M)\subseteq\V_A$. Since the core sections of $\V_A$ are all trivial, this shows that
the Lie bracket of $TA\to TM$ restricts to $\V_A\to \V_M$.
\end{example}

\begin{lemma}\label{char_par_sections_nabla_A}
Let $(A,F_A)$ be a foliated Lie algebroid and  $\nabla^A$ the connection 
defined in\eqref{nabla_A}. The following are equivalent for
$a\in \Gamma(A)$:
\begin{enumerate}
\item the class $\bar a$  is  $\nabla^{A}$-parallel,
\item    $Ta|_{F_M}$ is a section of $F_A$,
\item the core section $a^\uparrow$ of $F_A\to A$ is parallel with respect 
to the Bott connection $\nabla^{F_A}$, that is $\left[a^\uparrow,\Gamma(F_A)\right]\subseteq \Gamma(F_A)$.
\end{enumerate} 

\end{lemma}

\begin{proof}
The equivalence $(2)\Leftrightarrow (3)$ has been shown in Lemma \ref{lem_a_uparrow_preserves_FA}
and the implication $(3)\Rightarrow (1)$ is obvious.

If $\nabla^{A}_{\bar X}\bar a=0$ for all $\bar X\in\Gamma(F_M)$,
we get that 
\[  \left[X, a^\uparrow\right]\in\Gamma(F_A)\]
for all linear sections $X\in\Gamma(F_A)$. Since 
\[\left[a^\uparrow, b^\uparrow\right]=0\in\Gamma(F_A)
\]
for all $b\in\Gamma(F_\core)$ and 
the linear sections of $F_A$ together with the core sections of $F_A$ span $F_A$, we get
\[  \left[a^\uparrow, \Gamma(F_A)\right]\subseteq \Gamma(F_A)\]
and the flow of $a^\uparrow$ leaves 
$F_A$ invariant. 

Hence, we have shown $(1)\Rightarrow (3)$.
\end{proof}

Now we can complete the proof of the main theorem of this subsection.
\begin{theorem}\label{foliated_alg_yields_IM}
Let $(A,F_A)$ be a foliated Lie algebroid. Then the 
partial $F_M$-connection
$\nabla^{A}$ on $A/F_\core$ as in \eqref{nabla_A}  has the same properties as the connection $\nabla$  in 
Theorem \ref{the_connection}, and such that
\[ \bar a \text{ is } \nabla^{A}\text{-parallel}
\quad \Leftrightarrow \quad  Ta(F_M)\subseteq F_A.
\]
\end{theorem}

\begin{proof}
By Proposition \ref{nabla_A_welldef}, we know that 
$\nabla^A$ as in \eqref{nabla_A} is well-defined and flat and 
we have shown in Lemma \ref{char_par_sections_nabla_A} that 
$\bar a$ is $\nabla^{A}$-parallel if and only if 
$Ta(F_M)\subseteq F_A$.

If $c\in\Gamma(F_\core)$, then $\hat c\an{F_M}$ is a section of $F_A\to F_M$, 
and $\widehat{[a,c]}\an{F_M}=[Ta\an{F_M},\hat c\an{F_M}]\in\Gamma_{F_M}(F_A)$.
Hence, $[a,c]\in\Gamma(F_\core)$.

If $b\in\Gamma(A)$ is also such that $\bar b$ is $\nabla^A$-parallel, and 
$Ta\an{F_M}, Tb\an{F_M}\in\Gamma_{F_M}(F_A)$, then, since 
$F_A\to F_M$ is a subalgebroid of $TA\to TM$, we find 
that $[Ta\an{F_M}, Tb\an{F_M}]\in\Gamma_{F_M}(F_A)$.
But since $[Ta\an{F_M}, Tb\an{F_M}]=T[a,b]\an{F_M}$, this shows that
$\overline{[a,b]}$ is $\nabla^A$-parallel by Lemma \ref{char_par_sections_nabla_A}.

It remains to show that $\rho(a)\in\mx(M)$ is $\nabla^{F_M}$-parallel
if $a$ is $\nabla^A$-parallel.
Since $Ta\an{F_M}\in\Gamma_{F_M}(F_A)$
and $F_A\to F_M$ is a subalgebroid of $TA\to TM$, we find that
\[\rho_{TA}(T_ma(v_m))=J_M(T_m(\rho_A(a))v_m)\]
 is  an element of $T_{v_m}F_M$ for any $v_m\in F_M$.
Choose $X\in\Gamma(F_M)$ and  $\alpha\in\Gamma(F_M^\circ)$. Then, if $l_\alpha\in C^\infty(TM)$
is the linear function defined by $\alpha$,  we have 
\[\dr_{X(m)}l_\alpha(\rho_{TA}(T_ma(X(m))))=0.\] But the identity
\[J_M(T_m(\rho_A(a))X(m))=\left.\frac{d}{dt}\right\an{t=0}T_m\phi_tX(m),\] 
where $\phi_\cdot$ is the flow of $\rho_A(a)$,
 yields then immediately
\[\left(\ldr{\rho_A(a)}\alpha\right)(X(m))=0\]
for all $m$, and the equality
\[0=\ldr{\rho_A(a)}(\alpha(X))=(\ldr{\rho_A(\alpha)}\alpha)(X)+\alpha([\rho_A(a),X])
\]
leads thus to $\alpha([\rho_A(a),X])=0$.
Since $\alpha$ and $X$ were arbitrary sections of $F_M^\circ$ and $F_M$, respectively, 
we have shown that $[\rho_A(a),\Gamma(F_M)]\subseteq \Gamma(F_M)$.
\end{proof}

\begin{theorem}
Let $(A,F_A)$ be a foliated algebroid.
Assume that $A$ integrates to a Lie groupoid $G\rr M$ and that $j_G(F_A)$ integrates to $F_G\subseteq TG$.
Then 
\[ \nabla=\nabla^\s=\nabla^{A}.
\]
\end{theorem}

\begin{proof}
It is clear that if $j_G(F_A)$
integrates to $F_G$, then the core and the base of $F_G$ and 
of $F_A$ coincide. Hence, the three connections are flat partial $F_M$-connections 
on $A/F_\core$. Since $\nabla=\nabla^\s$
by Corollary \ref{char_of_nabla_s_parallel_sections_cor}, 
we only need to show that $\nabla$ and $\nabla^A$ have the same parallel sections.
For that we use Theorem \ref{char_of_nabla_s_parallel_sections}.
The class $\bar a$ of $a\in\Gamma(A)$ is $\nabla$-parallel
if and only if $\beta_a\an{F_M}$ is a section of $A(F_G)$.
Since $F_A=j_G\inv(A(F_G))$ and $j_G\inv\circ \beta_a\an{F_M}=Ta\an{F_M}$, 
this is equivalent to $Ta\an{F_M}\in\Gamma_{F_M}(F_A)$. But by the preceding theorem, this is true 
if and only if $\bar a$ is $\nabla^A$-parallel. 
\end{proof}

\vspace*{0.5cm}

\section{The connection as the infinitesimal data of a multiplicative foliation}
\label{connection-as-infinitesimal-data}
\subsection{Infinitesimal description of a multiplicative foliation}
We have seen that both foliated groupoids and foliated algebroids induce partial
 connections which have important properties. This motivates the following concept.
\begin{definition}\label{def_foliated_algebroid}
Let $(A\to M, \rho,[\cdot\,,\cdot])$ be a Lie algebroid, 
$F_M\subseteq TM$ an involutive subbundle,
$F_\core\subseteq A$ a subalgebroid over $M$ such that $\rho(F_\core)\subseteq F_M$
and $\nabla$ a partial $F_M$-connection on $A/F_\core$ with the following properties:
\begin{enumerate}
\item $\nabla$ is flat.
\item If $a\in\Gamma(A)$ is $\nabla$-parallel, then $[a,b]\in\Gamma(F_\core)$ 
for all $b\in\Gamma(F_\core)$.
\item If $a,b\in\Gamma(A)$ are $\nabla$-parallel, then $[a,b]$ is also parallel.
\item If $a\in\Gamma(A)$ is  $\nabla$-parallel, then $[\rho(a),\bar X]\in\Gamma(F_M)$ for all
$\bar X\in\Gamma(F_M)$. That is, $\rho(a)$ is $\nabla^{F_M}$-parallel.
\end{enumerate}
The quadruple $(A,F_M,F_\core,\nabla)$ will be referred to as an \textbf{IM-foliation\footnote{Here, ``IM'' stands for ``infinitesimal multiplicative'', as for the IM-$2$-forms in 
\cite{BuCr05}.}
 on $A$}, or simply an \textbf{IM-foliation}.
\end{definition}

\begin{remark}
\begin{enumerate}
\item IM-foliations already appear in the work of Eli Hawkins \cite{Hawkins08}, 
where they are also found to correspond with foliated algebroids.
\item Note that this is an infinitesimal version of the \textbf{ideal systems} in \cite{Mackenzie05}.
\end{enumerate}
\end{remark}

We have shown in Theorem \ref{the_connection} that each foliated Lie groupoid $(G\rr M, F_G)$
 induces an IM-foliation
on $A$. In the same manner, we have seen in Theorem \ref{foliated_alg_yields_IM}
that each foliated Lie algebroid $(A,F_A)$ induces an IM-foliation
on $A$. Moreover, 
a foliated Lie groupoid $(G\rr M, F_G)$ and its foliated Lie algebroid $(A,j_G\inv(A(F_G)))$
induce the same IM-foliation on $A$.

We will show here that an IM-foliation on a Lie algebroid is equivalent 
to a morphic foliation on this Lie algebroid.
\emph{This implies that IM-foliations on integrable Lie algebroids 
are in one-to-one correspondence with 
multiplicative foliations on the corresponding Lie groupoids.}

\subsection{Reconstruction of the foliated Lie algebroid from the connection}

Let now  $(q_A:A\to M,\rho, [\cdot\,,\cdot])$ be a Lie algebroid 
and $(A,F_M,F_\core,\nabla)$ an IM-foliation on $A$.

The Lie algebroid $TA\to TM$ is spanned by the sections 
$Ta$ and $\hat a$ for all $a\in\Gamma(A)$. 
Recall that
 $\hat a(v_m)=T_m0^Av_m+_{p_A}\left.\frac{d}{dt}\right\an{t=0}ta(m)$
for $v_m\in T_mM$ and that the bracket of the Lie algebroid 
$TA\to TM$ is given by
$$\left[Ta,Tb\right]_{TA}=T[a,b],\qquad 
\left[Ta,\hat b\right]_{TA}=\widehat{[a,b]},
\qquad \left[\hat a,\hat b\right]_{TA}=0$$
for all $a,b\in\Gamma(A)$.

For each $v_m\in F_M$, define the subset $F_A(v_m)$ of $(TA)_{v_m}$ as follows:
\begin{equation}\label{def_of_FA}
F_A(v_m)=\erz\left\{\begin{array}{c}
T_ma(v_m)\in T_{a(m)}A,\\[1em] 
\hat b(v_m)
\end{array}\left|\begin{array}{c}  
a\in\Gamma(A) \text{ is }\nabla\text{-parallel, }\\
b\in\Gamma(F_\core)
\end{array}
\right.\right\}.
\end{equation}
That is, $F_A$ is spanned by the restrictions to $F_M$ of the linear sections $Ta$ of $TA$ 
defined by $a\in\Gamma(A)$ such that $\bar{a}\in\Gamma(A/F_\core)$ is
 $\nabla$-parallel and the core sections $\hat b$ defined 
by sections $b\in\Gamma(F_\core)$.
We will show that $F_A$ is a well-defined  subalgebroid of $TA\to TM$.

 First, we check that 
$F_A$ is a double vector bundle
\begin{displaymath}\begin{xy}\label{diagram}
\xymatrix{ 
 F_A\ar[d]\ar[r]&A\ar[d]\\
F_M\ar[r]&M
}
\end{xy}.
\end{displaymath}
By construction, $F_A$ has core $F_\core$ and 
there will be 
 an injective double vector bundle morphism to
\begin{displaymath}\begin{xy}\label{diagram}
\xymatrix{ 
 TA\ar[d]_{Tq_A}\ar[r]^{p_A}&A\ar[d]^{q_A}\\
TM\ar[r]_{p_M}&M
}
\end{xy}
\end{displaymath}
which has core $A$.

By definition, we have $T_{a_m}q_A(T_ma(v_m))
=T_m(q_A\circ a)v_m=v_m\in F_M(m)$ for all $a\in\Gamma(A)$ and 
$v_m\in F_M(m)$ and 
$T_{v_m}q_A\left(\hat b(v_m)\right)=v_m\in F_M(m)$
for all $b\in\Gamma(F_\core)$.

Choose $m\in M$. Then there exists a neighborhood 
$U$ of $m$ in $M$ and 
sections $a_1, \ldots,a_n\in\Gamma(A)$
such that $\overline{a_{r+1}},\ldots,\overline{a_n}$ are $\nabla$-parallel and
form a basis for  $A/F_\core$ on $U$, and 
$a_{1},\ldots,a_r$ form a basis for $F_\core$ on $U$. Note that 
$\overline{a_{1}}, \ldots, \overline{a_r}$ are in a trivial manner $ \nabla$-parallel.
Choose any $\nabla$-parallel section $a\in\Gamma(A)$. Then, on $U$, 
we have $a=\sum_{i=1}^nf_ia_i$ 
with $f_{1},\ldots,f_r\in C^\infty(U)$ 
and $F_M$-invariant functions  $f_{r+1},\ldots,f_n\in C^\infty(U)$.
For any $v_m\in F_M(m)$, the equality 
\begin{align*}
T_ma(v_m)&=\sum_{i=1}^rv_m(f_i)\widehat{a_i}(v_m)+\sum_{i=1}^nf_i(m)T_ma_i(v_m)
\end{align*}
can be checked in coordinates, using $v_m(f_i)=0$ for $i=r+1,\ldots,n$.
If \[\sum_{i=1}^n\alpha_iT_ma_i(v_m)+\sum_{i=1}^r\beta_i\widehat{a_{i}}(v_m)=0^{TA}_{v_m}
\]
for some $v_m\in F_m$ and $\alpha_1,\ldots,\alpha_n,\beta_{1},\ldots,\beta_r\in\R$, 
we find 
\[\sum_{i=1}^n\alpha_ia_i(m)
=p_A\left(\sum_{i=1}^n\alpha_iT_ma_i(v_m)+\sum_{i=1}^r\beta_i\widehat{a_{i}}(v_m)\right)
=0^{A}_{m}.
\]
Since $a_1,\ldots,a_n$ is a local frame for $A$, this yields 
$\alpha_1=\ldots=\alpha_n=0$ and the equality $\beta_{1}=\ldots=\beta_r=0$ follows easily.
The set of sections $Ta_1,\ldots,Ta_n$
and $\widehat{a_{1}},\ldots, \widehat{a_r}$ is thus a local frame 
for $F_A$ as a vector bundle over $F_M$. 

\medskip

We next show that $F_A$
is also a vector bundle over $A$ (a subbundle of $TA\to A$).
Choose $a_m\in A_m$. Then there exists a $\nabla$-parallel
section $a\in\Gamma(A)$ defined on a neighborhood $U$ of $m$ such that $a(m)=a_m$.
Let $r$ be the rank of $F_\core$ and $l$ the rank of $F_M$.
Choose a basis of vector fields  $V_1,\ldots,V_l$ for $F_M$ on $U$ and a basis of sections
 $b_1,\ldots,b_{r}$ for $F_\core$ on $U$.  Choose a local frame $a_1,\ldots,a_n$ for 
$A$ on $U$ by $\nabla$-parallel sections 
and define for $i=1,\ldots,l$ the sections 
$\tilde V_i:q_A\inv(U)\to F_A$
by
\[\tilde V_i(a_m)=T\left(\phi^{a_1^\uparrow}_{\alpha_1}\circ\ldots\circ\phi^{a_l^\uparrow}_{\alpha_l}\circ 0^A \right)(V_i(m))
\]
where $\alpha_1,\ldots,\alpha_n\in\R$ are such that $a_m=\sum_{i=1}^n\alpha_i\cdot a_i(m)$. 
By construction, the vector fields $\tilde V_i$ are sections 
of $F_A\to A$ (recall Lemma \ref{char_par_sections_nabla_A}), 
and, using the fact that core vector fields commute, it is easy to show that
$\tilde V_i$ is linear over $V_i$ for each $i=1,\ldots,l$.
The sections $\tilde V_1,\ldots,\tilde V_l,b_1^\uparrow,\ldots,b_{r}^\uparrow$  
form  a basis for 
$F_A$ (seen as vector bundle over $A$) on $q_A\inv(U)$.

\medskip

Now we check that $(F_A,\rho_{TA},[\cdot\,,\cdot]_{TA})$
is a Lie algebroid over $F_M$ (a subalgebroid of $TA\to TM$).
Choose two linear sections $Ta\an{F_M}$ and $Tb\an{F_M}$ of 
$F_A$, i.e., with $\bar a, \bar b\in\Gamma(A/F_\core)$ that are $\nabla$-parallel.
We have then $[Ta,Tb]_{TA}=T[a,b]$. By the properties of 
the connection, $[a,b]$ is $\nabla$-parallel and $T[a,b]\an{F_M}$ is a section of $F_A$.
Choose now a core section $\hat b\an{F_M}$ and a linear section 
$Ta\an{F_M}$ of $F_A$, i.e., with $b\in\Gamma(F_\core)$ and 
$a\in\Gamma(A)$ a $\nabla$-parallel section.
We get 
$[Ta,\hat b]_{TA}=\widehat{[a,b]}$, whose restriction to $F_M$ is a section of $F_A$
since, by the properties of the connection, $[a,b]\in\Gamma(F_\core)$. 
Since the bracket of two core sections vanishes, we have shown that $[\cdot\,,\cdot]_{TA}$ 
restricts to $F_A$.

Next, we show that the anchor map of $TA$ restricts to a map 
$F_A\to TF_M$. Again, it is sufficient to show that the spanning sections 
of $F_A$ are sent by $\rho_{TA}$ to vector fields on $F_M$.
Recall that the anchor map is given by  
$\rho_{TA}=J_M\circ T\rho$. Choose a $\nabla$-parallel section 
$a\in\Gamma(A)$, set $X:=\rho(a)\in\mx(M)$  and compute
for any $v_m=\dot c(0)\in F_M$:
\begin{align*}
\rho_{TA}(T_ma(v_m))&=J_M(T_m(\rho(a))v_m)
=\left.\frac{d}{ds}\right\an{s=0}T_m\phi_s^X(v_m).
\end{align*}
By the properties of the connection, the vector field $X$
is such that $[X,\Gamma(F_M)]\subseteq \Gamma(F_M)$. 
By Corollary \ref{appendix_cor}, this implies $T_m\phi_s^XF_M(m)=F_M(\phi_s^X(m))$ for all $m\in M$ and 
$s\in\R$ where this makes sense.
 This implies that the curve 
$s\mapsto T_m\phi_s^X(v_m)$ has image in $F_M$ and consequently, 
$\rho_{TA}(T_ma(v_m))\in T_{v_m}F_M$.

Similarly, choose $b\in\Gamma(F_\core)$, $v_m=\dot c(0)\in F_M(m)$, set $Y=\rho(b)\in\Gamma(F_M)$
  and 
compute:
\begin{align*}
&\rho_{TA}\left(T_m0^A(v_m)+_{p_A}\left.\frac{d}{dt}\right\an{t=0}tb(m)\right)\\
=\,&
J_M\left(T_{0_m}\rho\left(T_m0^A(v_m)+_{p_A}\left.\frac{d}{dt}\right\an{t=0}tb(m)\right)\right)\\
=\,&J_M\left(\left.\frac{d}{dt}\right\an{t=0}0^{TM}(c(t))+_{p_{TM}}
\left.\frac{d}{dt}\right\an{t=0}tY(m)\right)\\
=\,&J_M\left(\left.\frac{d}{dt}\right\an{t=0}\left.\frac{d}{ds}\right\an{s=0}c(t)
+_{p_{TM}}\left.\frac{d}{dt}\right\an{t=0}\left.\frac{d}{ds}\right\an{s=0}\phi_{st}^Y(m)\right)\\
=\,&\left.\frac{d}{ds}\right\an{s=0}\left.\frac{d}{dt}\right\an{t=0}c(t)
+_{Tp_M}\left.\frac{d}{ds}\right\an{s=0}\left.\frac{d}{dt}\right\an{t=0}\phi_{st}^Y(m)\\
=\,&\left.\frac{d}{ds}\right\an{s=0}v_m
+_{TM}sY(m)\in T_{v_m}F_M.
\end{align*}

\medskip

Note also that the induced map 
\begin{align*}
\begin{xy}
\xymatrix{
F_A\ar[r]^{p_A}\ar[d]_{Tq_A}& A\ar[d]^{q_A}\\
F_M\ar[r]_{p_M}&M
}
\end{xy}
\end{align*} is 
a surjective Lie algebroid morphism by construction.

We have proved the following proposition.
\begin{proposition}
Let $(A\to M, F_M, F_\core,\nabla)$ be an IM-foliation on $A$. Consider $F_A\subseteq TA$ constructed 
as in \eqref{def_of_FA}. Then $F_A\to F_M$ is a subalgebroid of $TA\to TM$.
\end{proposition}

We show now that $F_A\to A$ is also a subalgebroid of $TA\to A$.
\begin{proposition}
Let $(A\to M, F_M, F_\core,\nabla)$ be an IM-foliation on $A$. Consider $F_A\to F_M$ the Lie subalgebroid of 
$TA\to TM$ as above.
The vector subbundle $F_A\to A$ of $TA\to A$ is involutive. 
\end{proposition}
\begin{proof}
By definition, the linear sections of $TA\to TM$ that restrict to sections 
of $F_A\to F_M$ are the sections $Ta$ for all $a\in\Gamma(A)$ such that
$\bar a$ is a  $\nabla$-parallel section of $A/F_\core$.
Note also that, by definition, every section of $F_\core$ is $\nabla$-parallel in this sense.

As in the proof of Lemma \ref{lem_a_uparrow_preserves_FA}, this implies
$[a^\uparrow,\Gamma(F_A)]\subseteq \Gamma(F_A)$  if $a\in\Gamma(A)$ is a $\nabla$-parallel section.
  As a consequence,  for any linear section
$X$   of $F_A$ over $\bar X$, we have 
$[X,a^\uparrow]=(D_Xa)^\uparrow$  with $D_Xa\in\Gamma(F_\core)$.
We find  thus the trivial equality $\nabla_{\bar X}\bar a=\overline{D_Xa}$ for  $\nabla$-parallel
sections $a\in\Gamma(A)$. 

Since $\nabla$ is flat, 
there exist local frames for $A$ of $\nabla$-parallel sections.
Choose $a\in\Gamma(A)$ and write $a=\sum_{i=1}^nf_ia_i$
with $a_1,\ldots,a_n\in\Gamma(A)$ $\nabla$-parallel sections and $f_1,\ldots,f_n\in C^\infty(M)$.
Then $\nabla_{\bar X}\bar a=\sum_{i=1}^n\bar X(f_i)\bar{a_i}$. 
But, on the other hand, if $X\in\Gamma(F_A)$ is a linear section over $\bar X$, we have 
\[(D_Xa)^\uparrow=\left[X,\left(\sum_{i=1}^nf_ia_i\right)^\uparrow\right]
=\sum_{i=1}^nq_A^*\left(\bar X(f_i)\right)a_i^\uparrow+\sum_{i=1}^nq_A^*f_i(D_Xa_i)^\uparrow
,\]
which leads to
\begin{align*}
\nabla_{\bar X}\bar a&=\sum_{i=1}^n\bar X(f_i)\bar{a_i}
=\sum_{i=1}^n\bar X(f_i)\bar{a_i}+\sum_{i=1}^nf_i\cdot \overline{D_Xa_i}=\overline{D_Xa},
\end{align*}
since $D_Xa_i\in\Gamma(F_\core)$ for $i=1,\ldots,n$.

Choose now two linear sections $(X,\bar X)$ and $(Y,\bar Y)$ of $F_A\to A$.
Since $\nabla$ is flat and $F_M$ is involutive , we have 
$\nabla_{[\bar X,\bar Y]}\bar a=\nabla_{\bar X}(\nabla_{\bar Y}\bar a)-\nabla_{\bar Y}(\nabla_{\bar X}\bar a)
$ for all $a\in\Gamma(A)$.
Hence, if
 $\nabla_{[\bar X,\bar Y]}\bar a=\bar k$ for some
$k\in\Gamma(A)$,  
then 
\[(k+b)^\uparrow=\left[X,\left[Y,a^\uparrow\right]\right]-\left[Y,\left[X,a^\uparrow\right]\right]
=\left[[X,Y],a^\uparrow\right]\]
for some $b\in\Gamma(F_\core)$
and, if $Z\sim_{q_A}[\bar X,\bar Y]$
is a linear section of $F_A$ covering $[\bar X,\bar Y]$, then 
\[(k+b')^\uparrow=\left[Z,a^\uparrow\right]\]
for some $b'\in\Gamma(F_\core)$. The linear section $W:=Z-[X,Y]$ covers hence $0^{TM}$
and satisfies 
\[ \left[W,a^\uparrow\right]=(b-b')^\uparrow.\]
 Since $b-b'\in\Gamma(F_\core)$, we have found
$D_Wa\in\Gamma(F_\core)$ for any $a\in\Gamma(A)$, and,
 combining this with Lemma 
\ref{linear_over_0}, one concludes
that $W\in\Gamma(F_A)$ and hence $[X,Y]\in\Gamma(F_A)$.
Since $F_A\to A$ is spanned by its linear and core sections, we have shown 
that $F_A$ is involutive. 
\end{proof}

This completes the proof of our  main  theorem.
\begin{theorem}\label{main_theorem}
Let $(q_A:A\to M,\rho_A,[\cdot\,,\cdot]_A)$ be a Lie algebroid
and $(A,F_M,F_\core,\nabla)$ an IM-foliation on $A$.
Then $F_A$ defined as in \eqref{def_of_FA} is a $\mathcal{VB}$-algebroid 
with sides $F_M$ and $A$ and core $F_\core$, such that
the inclusion 
is a $\mathcal{VB}$-algebroid morphism:
\begin{displaymath}\begin{xy}
\xymatrix{ 
F_A\ar@{^{(}->}[dr]\ar[rr]\ar[dd]&&A\ar@2{-}[dr]\ar[dd]&\\
&TA\ar[dd]\ar[rr]&&A\ar[dd]\\
F_M\ar@{^{(}->}[dr]\ar[rr]&&M\ar@2{-}[dr]&\\
&TM\ar[rr]&&M
}
\end{xy}
\end{displaymath}
Hence, $(A,F_A)$ is a foliated Lie algebroid and 
 the connection $\nabla^{A}$ induced by $F_A$ is equal to $\nabla$.
\end{theorem}

\begin{example} 
Assume that $H$ acts on a Lie groupoid $G$ over $M$ by groupoid automorphisms. 
Assume also that the action is free and proper. Starting from the data $(A,\mathcal{V}_M,0,\nabla)$ 
where $\nabla$ is the partial $\mathcal{V}_M$-connection on $A$ determined by
\[[\xi_G,a^r]=(\nabla_{\xi_M}a)^r\]
for all $\xi\in\lie g$ and $a\in\Gamma(A)$,
  the last theorem states that we  recover exactly the Lie foliated algebroid $\mathcal{V}_A\rmap \mathcal{V}_M$ 
obtained by applying the Lie functor to the foliated groupoid $\mathcal{V}_G\rightrightarrows \mathcal{V}_M$.
\end{example}

\begin{example} 
Assume that $\mathfrak{g}$ is a Lie algebra, i.e. a Lie algebroid over a point. 
In this case, the tangent Lie algebroid $T\mathfrak{g}$ is also a Lie algebroid over a point, 
that is, $T\mathfrak{g}$ is a Lie algebra. It is easy to see
 that the Lie algebra structure on $T\mathfrak{g}=\mathfrak{g}\times \mathfrak{g}$
 is the semi-direct 
product Lie algebra $\mathfrak{g}\ltimes \mathfrak{g}$ with respect to the adjoint representation 
of $\mathfrak{g}$ on itself. Note also that the fact that a quadruple $(\mathfrak{g},0,\mathfrak{f},\nabla=0)$ 
is an IM-foliation on $\mathfrak{g}$ is equivalent to saying that $\mathfrak{f}\subseteq \mathfrak{g}$ 
is an ideal. 

The morphic foliation $F_{\mathfrak{g}}$ associated to 
the IM-foliation $(\mathfrak{g},0,\mathfrak{f},\nabla=0)$ is given by 
$F_{\mathfrak{g}}=\mathfrak{g}\times \mathfrak{f}$. This follows immediately from the fact 
that every element $a\in \mathfrak{g}$ can be viewed, in a trivial way, as a $\nabla$-parallel 
section of $\mathfrak{g}\to \{0\}$. The property that $F_{\mathfrak{g}}=\mathfrak{g}\times \mathfrak{f}$ 
is a morphic foliation is equivalent to saying that $\mathfrak{g}\times \mathfrak{f}$ is a Lie 
subalgebra of $\mathfrak{g}\ltimes \mathfrak{g}$. In particular, if $G$ is the connected and 
simply connected Lie group integrating $\mathfrak{g}$, we conclude that the foliated algebroid 
$F_{\mathfrak{g}}=\mathfrak{g}\times \mathfrak{f}$ integrates to a Lie subgroup $G\times \mathfrak{f}$ 
of the semi-direct Lie group $G\ltimes \mathfrak{g}$ determined by the adjoint action of $G$ on its 
Lie algebra $\mathfrak{g}$. Using right (or left) translations, we get a subbundle $F_G\subseteq TG$ 
which is involutive and multiplicative. We conclude that $F_G$ is the multiplicative foliation on $G$ 
integrating the IM-foliation $(\mathfrak{g},0,\mathfrak{f},\nabla=0)$ on $\mathfrak{g}$. Thus, in the 
case of Lie groups and Lie algebras, this recovers the results in
\cite{Ortiz08, Jotz11a,Jotz11b}.
\end{example}

\vspace*{0.5cm}

\section{The leaf space of a foliated algebroid}\label{quotient_algebroid}
Assume that $(A, F_M, F_\core,\nabla)$ is an IM-foliation on $A$.
Then there is an induced involutive subbundle $F_A\subseteq TA$ as in Theorem 
\ref{main_theorem}.
We will show that if the leaf space $M/F_M$ is a  smooth manifold, 
such that the quotient map
$\pi_{M}:M\to M/F_M$ is a  
surjective submersion, then there is an induced Lie algebroid structure
$([q_A]: A/F_A\to M/F_M, [\rho],[\cdot\,,\cdot]_{A/F_A})$ such that the projection 
\begin{align*}
\begin{xy}
\xymatrix{
A\ar[d]_{q_A}\ar[r]^{\pi}&A/F_A\ar[d]^{[q_A]}\\
M\ar[r]_{\pi_M}&M/F_M
}
\end{xy}
\end{align*}
is a Lie algebroid morphism. Furthermore, if $A\to M$ integrates to a Lie groupoid $G\rr M$
and the completeness and regularity conditions for the leaf space $G/F_G\rr M/F_M$
to be a Lie groupoid are satisfied (see \cite{Jotz11b}), 
then $A/F_A\to M/F_M$ is the Lie algebroid of $G/F_G\rr M/F_M$.
We will see in the proofs that the important data for this reduction process is the IM-foliation $(A,F_M,F_\core,\nabla)$.

The class of $a_m\in A_m$ is written  $[a_m]\in A/F_A$, and in the same manner,
 the class of $m\in M$ is denoted by 
$[m]\in M/F_M$.
The class of $ a_m\in A_m$ in $A/F_\core$ will be written $\bar a_m$.
\begin{proposition}\label{prop_structure_maps_1}
Let $(A\to M, F_M, F_\core,\nabla)$ be an IM-foliation on $A$ and $F_A\subseteq TA$ 
the corresponding morphic foliation as in Theorem 
\ref{main_theorem}.
\begin{enumerate}
\item The map
$\pi:A\to A/F_A$ factors as a composition
\begin{align*}
\begin{xy}
\xymatrix{
A\ar[d]\ar[dr]^{\pi}&\\
A/F_\core\ar[r]_{\pi_\core}&A/F_A
}
\end{xy}
\end{align*}
That is, we have $\pi(a_m+b_m)=\pi(a_m)$
for all $a_m\in A$ and $b_m\in F_\core(m)$.
\item The equivalence relation $\sim:=\sim_{F_A}$ on $A$ can be described as follows.
\begin{equation}
a_m\sim a_n
\Leftrightarrow 
\begin{array}{c}
\text{ There exist linear sections }
(X_1,\bar X_1), \ldots, (X_r,\bar X_r) \text{ of } F_A\to A\\
\text{ with flows } \phi^1,\ldots, \phi^r
\text{ such that }\\
a_m\in \phi_{t_1}^1\circ\ldots\circ\phi_{t_r}^r(a_n)+F_\core(m)\\
\text{ for  some } t_1,\ldots, t_r\in \R.
\end{array}
\end{equation}

\item The map $q_A$ induces a map $[q_A]:A/{F_A}\to M/F_M$
such that 
\begin{align*}
\begin{xy}
\xymatrix{
A\ar[d]_{q_A}\ar[r]^{\pi}&A/{F_A}\ar[d]^{[q_A]}\\
M\ar[r]_{\pi_{M}}&M/F_M
}
\end{xy}
\end{align*}
commutes.
\item Let $a\in\Gamma(A)$ be such that $\bar a\in\Gamma(A/F_\core)$
is 
$\nabla$-parallel. 
Let $U:=\dom(a)\subseteq M$. Then there is an induced map $[a]: \pi_{M}(U)\to A/F_A$
such that 
\begin{align*}
\begin{xy}
\xymatrix{
U\ar[d]_{\pi_{M}}\ar[r]^a&A\ar[d]^{\pi}\\
\pi_{M}(U)\ar[r]_{[a]}&A/{F_A}
}
\end{xy}
\end{align*}
commutes. The map $[a]$ is a section of $[q_A]$ in the sense that
\[[q_A]\circ[a]=\id_{\pi_{M}(U)}.\]
\end{enumerate}
\end{proposition}
\begin{proof}
\begin{enumerate}
\item  Recall that all the core sections $b^\uparrow\in\mx(A)$
with $b\in\Gamma(F_\core)$ are sections of $F_A$.
Choose $a_m\in A$ and $b_m\in F_\core(m)$. Then there exists a section 
$b\in\Gamma(F_\core)$ with $b(m)=b_m$. The flow $\phi^{b^\uparrow}$ of $b^\uparrow$ is given
by $\phi^{b^\uparrow}_t(a)=a+tb(q_A(a))$
for all $a\in A$ and $t\in\R$. Hence, we have $a_m\sim a_m+tb(m)=a_m+tb_m$
for all $t\in\R$, and in particular $a_m\sim a_m+b_m$.
The map $\pi_\core:A/F_\core\to A/F_A$, 
$\pi_\core(a_m+F_\core(m))=[a_m]$ is hence well-defined and the diagram commutes.
\item  Since the family
of linear sections of $F_A$ and the family of core sections of $F_A$ span together $F_A$,
its leaves  are the accessible sets of these two families of vector fields
 (\cite{Stefan74a, Stefan80, Sussmann73}, see \cite{OrRa04}
for a review of these results).  Hence, two
points $a_m$ and $a_n$ in $A$
are in the same leaf of $F_A$ if they can be joined by finitely many 
curves along flow lines of core sections $b^\uparrow$ for 
$b\in\Gamma(F_\core)$ and linear vector fields $X\in\Gamma(F_A)$. 
By the involutivity of $F_A$, we have 
$D_Xb\in\Gamma(F_\core)$ for all $b\in\Gamma(F_\core)$ and linear vector fields $X\in\Gamma(F_A)$.
Hence, by Lemma \ref{appendix_util}, we get that $F_\core$ is invariant under the flow lines of linear
vector fields with values in $F_A$.
That is, using the fact that $\phi^X_t$ is a vector bundle homomorphism, we have 
\[\left(\phi^X_t\circ\phi_s^{b^\uparrow}\right)(a_m)\in \phi^X_t(a_m+F_\core(m))
=\phi^X_t(a_m)+F_\core\left(\phi_t^{\bar X}(m)\right)
\]
for all $a_m\in A$, $t\in\R$ where this makes sense and $s\in\R$. Since 
\[\phi_s^{b^\uparrow}\circ\phi^X_t(a_m)\in \phi^X_t(a_m)+F_\core\left(\phi_t^{\bar X}(m)\right),
\]
 the proof is finished.
\item  Assume that $a_m\sim a_n$ for some elements $a_m,a_n\in A$. Then  there exists, 
without loss of generality, one linear vector field 
$X\in\Gamma(F_A)$ over $\bar X\in\Gamma(F_M)$, an element $b_m\in F_\core(m)$ and $t\in \R$
such that $a_m=\phi^X_t(a_n)+b_m$.
We have then immediately 
\[m=q_A(a_m)=q_A\left(\phi^X_t(a_n)+b_m\right)=\left(q_A\circ \phi^X_t\right)(a_n)=\phi^{\bar X}_t(n),
\]
which shows $m\sim_{F_M}n$.
\item Assume first that $a$ does not vanish on its domain of definition.
 Since $a$ is $\nabla$-parallel,  we have 
$D_Xa\in\Gamma(F_\core)$ for all linear vector fields $X\in\Gamma(F_A)$, $X\sim_{q_A}\bar X\in\Gamma(F_M)$.
By Lemma \ref{appendix_util}, this yields 
\begin{equation}\label{eq_for_a}
\phi^X_t(a(m))\in a\left(\phi^{\bar X}_t(m)\right)+F_\core\left(\phi^{\bar X}_t(m)\right)
\end{equation}
for all $t\in\R$ where this makes sense and consequently
\[\pi\left(a(m)\right)=\pi\left(a\left(\phi^{\bar X}_t(m)\right)\right).
\]
Since $F_M$ is spanned by projections $\bar X$
of  linear vector fields $X\in\Gamma(F_A)$, this shows that
$a$ projects to the map $[a]$.

In general, we have $a=\sum_{i=1}^nf_ia_i$ on an open set $U$ with non-vanishing 
$\nabla$-parallel sections $a_1,\ldots,a_n$ of $A$
such that $a_1,\ldots,a_r\in\Gamma(F_\core)$ 
and  functions $f_1,\ldots,f_n\in C^\infty(U)$ such that $f_{r+1},\ldots,f_n$ are 
$F_M$-invariant.
This yields using \eqref{eq_for_a}:
\begin{align*}
\phi^X_t(a(m))&=\phi^X_t\left(\sum_{i=1}^nf_i(m)a_i(m)
\right)=\sum_{i=1}^nf_i(m)\phi^X_t(a_i(m))\\
&\in\sum_{i=r+1}^nf_i\left(\phi^{\bar X}_t(m)\right)\phi^X_t(a_i(m))+F_\core\left(\phi^{\bar X}_t(m)
\right)\\
&\hspace*{4cm}=
a\left(\phi^{\bar X}_t(m)\right)+F_\core\left(\phi^{\bar X}_t(m)\right)
\end{align*}
and we get the statement in the same manner as above.

We have 
\[\left([q_A]\circ[a]\right)\circ\pi_M=
[q_A]\circ \pi\circ a
=\pi_M \circ q_A\circ a=\pi_M\circ \id_M=\pi_M,\]
which shows the last claim since $\pi_M$ is surjective.
\end{enumerate}
\end{proof}

\begin{corollary}\label{holonomy_cor}
Let $(A\to M, F_M, F_\core,\nabla)$ be an IM-foliation on a Lie
algebroid $A$.  Choose $\bar a_m$ and $\bar a_n$ in $A/F_\core$.
\begin{enumerate}
\item Then
$\pi_\core(\bar a_m)=\pi_\core(\bar a_n)$
if and only if $\overline{a_m}\in A/F_\core$ is the $\nabla$-parallel transport of $\overline{a_n}$ over 
a piecewise smooth path along the foliation defined by $F_M$ on $M$.
\item If $\nabla$ has trivial holonomy, then $\pi_\core(\overline{a_m})=\pi_\core(\overline{a_m'})$
if and only if $\overline{a_m}=\overline{a_m'}$.
\end{enumerate}
\end{corollary}
\begin{proof}
\begin{enumerate}
\item Assume first that $\pi_\core(\bar a_m)=\pi_\core(\bar a_n)$. Then there exists without loss of generality one 
 linear vector field $X\in\Gamma(F_A)$ over $\bar X\in\Gamma(F_M)$ and $t\in\R$
such that $\bar a_m=\overline{\phi_t^X(a_n)}$. 

Consider the curve $a:[0,t]\to A/F_\core$ over 
$c:=\phi_\cdot^{\bar X}(n)$ defined by
\[a(\tau)=\overline{\phi^X_\tau(a_n)}\]
for $\tau\in[0,t]$.
 For each $\tau\in[0,t]$, we find 
$\varepsilon_\tau>0$ and a parallel section $a^\tau$ of $A$ such that
$\overline{a^\tau(c(\tau))}=a(\tau)$.
As in  the proof of Proposition \ref{prop_structure_maps_1}, 4), we get 
then that $\overline{a^\tau(c(s))}=a(s)$ for $s\in[\tau-\varepsilon_\tau,\tau+\varepsilon_\tau]$.
This yields $\nabla_{\bar X(c(\tau))}a=0$ for 
all $\tau$.

\medskip

Conversely, assume that $\overline{a_m}\in A/F_\core$ is the $\nabla$-parallel transport of $\overline{a_n}$ over 
a piecewise smooth path along a path lying in the leaf of $F_M$
through $n$. Without loss of generality, this path is a segment of a
flow curve of a vector field $\bar X\in\Gamma(F_M)$, $m=\phi_t^{\bar
  X}(n)$ for some $t\in\R$, and 
there exists a $\nabla$-parallel section $a$ of $A$ such that
$\overline{a(m)}=\overline{a_m}$ and $\overline{a(n)}=\overline{a_n}$.
Choose any linear vector field $X\in\Gamma(F_A)$ over $\bar X$. Then 
we get as in the proof  of Proposition \ref{prop_structure_maps_1}, 4)
that 
\[\overline{a(m)}=\overline{a\left(\phi_t^{\bar X}(n)\right)}
=\overline{\phi_t^{X}(a(n))}\]
and hence $a_m\sim a_n$ by Proposition \ref{prop_structure_maps_1}, 2).
\item This is immediate since here, parallel transport does not depend
  on the path along the leaf of $F_M$ through $m$.
\end{enumerate}
\end{proof}

Using Proposition \ref{prop_structure_maps_1}, we show that 
if  $M/F_M$ is a smooth manifold and the projection is a  
submersion, then the quotient $A/F_A$ 
is  a vector bundle over $M/F_M$.

\begin{proposition}
Let $(A\to M, F_M, F_\core,\nabla)$ be an IM-foliation on a Lie
algebroid $A$. 
Assume that  $M/F_M$ is a smooth manifold such that the projection is a  
submersion, and that the connection $\nabla$ has no holonomy.
The quotient space $A/F_A$ inherits a vector bundle structure over $\bar M$
such that the projection $(\pi,\pi_M)$
is a vector bundle homomorphism.
\begin{displaymath}\begin{xy}
\xymatrix{ 
 A\ar[r]^{\pi}\ar[d]_{q_A}& A/F_A\ar[d]^{[q_A]}\\
M\ar[r]_{\pi_M}&M/F_M
}
\end{xy}
\end{displaymath}
\end{proposition}
\begin{proof}
Choose a local frame for $A/F_\core$ of $\nabla$-parallel sections
$\bar a_1,\ldots,\bar a_k$ defined on $U\subseteq M$, $k=\operatorname{rank}(A/F_\core)$.
Write $q_{A/F_\core}:A/F_\core\to M$ for the vector bundle projection and consider the local  trivialization
of $A/F_\core$:
\[\begin{array}{cccc}
\Phi:&q_{A/F_\core}\inv(U)&\to &U\times\R^k\\
& \overline{b_m}&\mapsto &(m, \xi_1(b_m),\ldots,\xi_k(b_m)),
\end{array}
\]
where $b_m=\sum_{i=1}^k\xi_i(b_m)a_i(m)+c_m$
with some $c_m\in F_\core(m)$.

By Corollary \ref{holonomy_cor}, we find that for
$\bar b_m$, $\bar b_n\in q_{A/F_\core}\inv(U)$, the equality 
\[\pi_\core\left(\overline{b_m}\right)=\pi_\core\left(\overline{b_n}\right)\]
implies 
\[\xi_i(b_m)=\xi_i(b_n) \text{ for } i=1,\ldots,k,\]
since $\bar b_m$ is the parallel transport of $\bar b_n$ 
along any path in the leaf of $F_M$ through $m$ and $n$, and so in particular
along a path in $U$.  
That is, the map $\Phi$ factors to a well-defined 
map 
\[[\Phi]:[q_A]\inv(\bar U)\to \bar U\times\R^k\]
such that 
\[[\Phi]\circ\pi_\core=(\pi_M,\id_{\R^k})\circ\Phi.\]

The map $[\Phi]$ is the projection to $A/F_A$ of the ``well chosen'' local trivialization
$q_{A/F_\core}\times\xi_1\times\ldots\times\xi_k$ of 
$A/F_\core$. Since by Proposition \ref{parallel_sections},
we can cover $A/F_\core$  by this type of $F_A$-invariant trivializations, 
we find that we can construct trivializations for $A/F_A$, which is hence 
shown to be a vector bundle over $\bar M$.
\end{proof}

\begin{remark}
Corollary \ref{holonomy_cor} implies that the 
quotient space $\pi_\core:A/F_\core\to A/F_A$
is the quotient by the equivalence relation
given by parallel transport. Constructions like this were made in 
\cite{Zambon08}, see also \cite{JoRaZa11}.
This idea will be used (implicitly) in the proofs of the following statements.

Note that this shows also that the data $(A,F_M,F_\core,\nabla)$
is an infinitesimal version of the ideal systems as in \cite{Mackenzie05},
and the methods of construction of the quotient algebroid 
are similar in this sense.
\end{remark}

By construction, we have also the following vector bundle homomorphism
\begin{displaymath}\begin{xy}
\xymatrix{ 
 A/F_\core\ar[r]^{\pi_\core}\ar[d]_{q_{A/F_\core}}&A/F_A\ar[d]^{[q_A]}\\
M\ar[r]_{\pi_M}&M/F_M
}
\end{xy}
\end{displaymath}
which is an isomorphism in every fiber,
and we find that for each local section $[a]$ of $A/F_A$ defined on $\bar U\subseteq \bar M$,
 there exists 
a $\nabla$-parallel section $a$ of $A$
defined on $\pi_M\inv(\bar U)$ such that $\pi\circ a=[a]\circ\pi_M$.

\begin{proposition}\label{the_anchor}
Let $(A,F_M,F_\core,\nabla)$ be an IM-foliation and assume 
that the quotient space $\bar M=M/F_M$
is a smooth manifold and $\nabla$ has trivial holonomy.
Then there is an induced map 
$[\rho]:A/F_A\to T\bar M$
such that 
\begin{displaymath}\begin{xy}
\xymatrix{ 
A\ar[r]^{\rho}\ar[d]_{\pi}&TM\ar[d]^{T\pi_M}\\
A/F_A\ar[r]_{[\rho]}&T\bar M
}
\end{xy}
\end{displaymath}
commutes.

If $a\in\Gamma(A)$
is $\nabla$-parallel, then 
$\rho(a)\in\mx(M)$ is $\nabla^{F_M}$-parallel
and  $\pi_M$-related to 
$[\rho][a]\in\mx(\bar M)$.
\end{proposition}
\begin{proof}
Define $[\rho]:A/F_A\to T(M/F_M)$ 
by 
\[[\rho]([a_m])=T_m\pi_M(\rho(a_m))\in T_{[m]}(M/F_M)\simeq T_mM/F_M(m).
\]
To see that $[\rho]$  is well-defined, recall first that 
$\rho(F_\core)\subseteq F_M$.
If $[a_m]=[a_n]$, then
$a_m=b_m+\phi^X_t(a_n)$
for some linear section $X\in\Gamma(F_A)$ over $\bar X\in\Gamma(F_M)$, $t\in\R$ 
and $b_m\in F_\core(m)$. 
Consider the curve $a:[0,t]\to A/F_\core$ over 
$c:=\phi_\cdot^{\bar X}(n)$ defined by
\[a(\tau)=\overline{\phi^X_\tau(a_n)}\]
for $\tau\in[0,t]$.
Then $\nabla_{\bar X(c(\tau))}a=0$ for 
all $\tau$. For each $\tau\in[0,t]$, we find 
$\varepsilon_\tau>0$ and a parallel section $a^\tau$ of $A$ such that
$a^\tau(c(s))=a(s)$ for $s\in[\tau-\varepsilon_\tau,\tau+\varepsilon_\tau]$.
Then, $\rho\circ a^\tau$ is $\nabla^{F_M}$-parallel, and we get by Lemma \ref{immediate_lemma}
that $\left(\rho\circ a^\tau\right)\sim_{\pi_M}Y^\tau$
for some $Y^\tau\in\mx(\bar M)$.

Since $[c(\tau-\varepsilon_\tau)]=[c(s)]$ for all $s\in[\tau-\varepsilon_\tau,\tau+\varepsilon_\tau]$, we have then 
\begin{align*}
T_{c(\tau-\varepsilon_\tau)}\pi_M(\rho(a(\tau-\varepsilon_\tau)))
&=T_{c(\tau-\varepsilon_\tau)}\pi_M(\rho(a^\tau(c(\tau-\varepsilon_\tau))))\\
&=Y^\tau([c(\tau-\varepsilon_\tau)])=Y^\tau([c(s)])
=T_{c(s)}\pi_M(\rho(a(s)))
\end{align*}
for all $s\in(\tau-\varepsilon_\tau,\tau+\varepsilon_\tau)$. Since $[0,t]$ is covered by intervals like this, we get
 $T_m\pi_M(\rho(a_m))=T_{m}\pi_M(\rho(a(0)))
=T_{n}\pi_M(\rho(a_n))$, which shows that $[\rho]$ is well-defined.
\end{proof}

Now we will define a Lie bracket on the space of sections of $A/F_A$. For $[a],[b]\in\Gamma(A/F_A)$, choose 
$\nabla$-parallel sections $a,b\in\Gamma(A)$
such that $[a]\circ\pi_M=\pi\circ a$ and
$[b]\circ\pi_M=\pi\circ b$. Then $[a,b]$ is $\nabla$-parallel
by the properties of $\nabla$ and
we can define 
\[\Bigl[[a],[b]\Bigr]_{A/F_A}\in\Gamma(A/F_A)\]
by
\[\Bigl[[a],[b]\Bigr]_{A/F_A}\circ \pi_M=\pi\circ[a,b].
\]
By the properties of $\nabla$, this definition does not depend on the choice of the $\nabla$-parallel
sections (which can be made up to sections of $F_\core$), 
and by definition and with Proposition \ref{the_anchor}, we get
\[[\rho]\left(\Bigl[[a],[b]\Bigr]_{A/F_A}\right)
=\Bigl[[\rho][a], [\rho][b]\Bigr]_{T\bar M},
\]
where the bracket on the right-hand side is the Lie bracket on vector fields on $\bar M$.

We can now complete the proof of the  following theorem.
\begin{theorem}
Let $(A,F_M,F_\core,\nabla)$ be an IM-foliation on a Lie algebroid
$A$. Assume that $\bar M=M/F_M$ is a smooth manifold and $\nabla$ has trivial holonomy.
Then the triple\linebreak $(A/F_A,[\rho],[\cdot\,,\cdot]_{A/F_A})$
is a Lie algebroid over $\bar M$
such that the projection 
$(\pi,\pi_M)$ is a Lie algebroid morphism.
\begin{align*}
\begin{xy}
\xymatrix{
A\ar[d]_{q_A}\ar[r]^{\pi}&A/F_A\ar[d]^{[q_A]}\\
M\ar[r]_{\pi}&M/F_M
}
\end{xy}
\end{align*}
\end{theorem}

\begin{proof}
The Jacobi identity follows immediately from the properties of the Lie bracket $[\cdot\,,\cdot]$ on $\Gamma_M(A)$ 
and the definition of $[\cdot\,,\cdot]_{A/F_A}$.
For the Leibniz identity  choose $[a],[b]\in\Gamma(A/F_A)$
corresponding to $\nabla$-parallel sections $a,b\in\Gamma_M(A)$,
and $f\in C^\infty(\bar M)$. We have then $\pi_M^*f\in C^\infty(M)^{F_M}$
and $f\cdot[b]$
corresponds  to the $\nabla$-parallel section
$(\pi_M^*f)\cdot b$ of $A$ (see Lemma \ref{immediate_lemma}).
We have hence 
\begin{align*}
[[a], f\cdot [b]]_{A/F_A}\circ \pi_M
&=\pi\circ [a, (\pi_M^*f)\cdot b]\\
&=\pi\circ \left(\pi_M^*f\cdot [a,b]+\rho(a)(\pi_M^*f)\cdot b\right)\\
&=\pi\circ \left(\pi_M^*f\cdot [a,b]+\pi_M^*([\rho][a](f))\cdot b\right)\\
&=\left(f\cdot [[a],  [b]]_{A/F_A}+[\rho][a](f)\cdot[b]\right)\circ \pi_M,
\end{align*}
where we have used Proposition \ref{the_anchor}
in the third equality.

The fact that $(\pi,\pi_M)$ is compatible with the Lie algebroid brackets 
is immediate  by construction
and the compatibility of the anchor maps is given by the definition of $[\rho]$.
\end{proof}

Assume now that $F_G\subseteq TG$ is a multiplicative foliation on $G\rr M$ 
such that the leaf space $G/F_G$ is a 
Lie groupoid over the leaf space $M/F_M$ (recall that there are topological conditions for this to be true 
\cite{Jotz11b}). The multiplicative foliation $F_G$ determines an IM-foliation 
$(A,F_M,F_\core,\nabla)$ and a Lie algebroid 
$(A/F_A,[\rho], [\cdot,\cdot]_{A/F_A})$ as in the preceding theorem. We conclude this subsection 
with the comparison of this Lie algebroid with the Lie algebroid of the leaf space $G/F_G\rr M/F_M$.

\begin{theorem}\label{lie_alg_lie_gpd}
Let $(A,F_M,F_\core,\nabla)$ be an IM-foliation on $A$.  Assume 
that $A$ integrates to a Lie groupoid $G\rr M$, and 
$F_A$  to a multiplicative foliation $F_G$ on $G$.
If $G/F_G$ and $M/F_M$ are smooth manifold, $\nabla$ has trivial holonomy and $F_G$ is such that 
$G/F_G\rr M/F_M$ is a Lie groupoid, then we have 
\[A(G/F_G)=A/F_A,
\]
where $A/F_A $ is equipped with the Lie algebroid structure in the previous theorem.
\end{theorem}

\begin{remark}
It would be interesting to study the relation between the trivial holonomy property of $\nabla$ and 
the condition of $F_G$ for $G/F_G\rr M/F_M$ to be a Lie groupoid.
\end{remark}

\begin{proof}[Proof of Theorem \ref{lie_alg_lie_gpd}]
Let $\pi_G:G\to G/F_G$ be the projection, and $[\s],[\tg]$
the source and target maps of $G/F_G\rr M/F_M$.
Recall from Theorem \ref{the_connection} that a section $a\in\Gamma(A)$
is $\nabla$-parallel if and only if 
$[a^r,\Gamma(F_G)]\subseteq \Gamma(F_G)$ and  the vector field $a^r$  is then $\pi_G$-related 
to a vector field $\overline{a^r}\in\mx(G/F_G)$.
We have 
\[T[\s]\circ \overline{a^r}\circ \pi_G=T[\s]\circ T\pi_G\circ a^r
=T\pi_G(T\s\circ a^r)=0,\]
which shows that $\overline{a^r}$ is tangent to the $[\s]$-fibers.
By Lemma 3.18 in \cite{Jotz11b}, 
we get 
\[\overline{a^r}([g])=T_g\pi_G(a^r(g))=T_g\pi_G(a(\tg(g))\star 0_g)
=T_{\tg(g)}\pi_G(a(\tg(g)))\star 0_{[g]}
=\overline{a^r}([\tg][g])\star 0_{[g]},
\]
which shows that $\overline{a^r}=\tilde a^r$
for  $\tilde a:=\overline{a^r}\an{M/F_M}\in\Gamma(A(G/F_G))$.

Since $(\pi_G,\pi_M)$ is a Lie groupoid morphism
\begin{align*}
\begin{xy}
\xymatrix{
G\ar[r]^{\pi_G}\ar@<.6ex>^{\s}[d]\ar@<-.6ex>_{\tg}[d]& G/F_G\ar@<.6ex>^{[\s]}[d]\ar@<-.6ex>_{[\tg]}[d] \\
M\ar[r]_{\pi_M}&M/F_M}
\end{xy}
\end{align*} 
the map
$A(\pi_G)=T\pi_G\an{A}$
\begin{align*}
\begin{xy}
\xymatrix{
A\ar[r]^{A(\pi_G)\,\,\,\,\,\,}\ar[d]& A(G/F_G)\ar[d] \\
M\ar[r]_{\pi_M}&M/F_M}
\end{xy}
\end{align*}
is a Lie algebroid morphism
and $A(\pi_G)(b_m)=0$
for all $b_m\in F_\core$.
For any  $\nabla$-parallel section
$a\in\Gamma(A)$  we have $a^r\sim_{\pi_G}\tilde a^r$ and hence 
\begin{equation}\label{the_equation}
A(\pi_G)\circ a =T\pi_G a=\tilde a\circ\pi_M.
\end{equation}
Define the map
\begin{align*}
\begin{xy}
\xymatrix{
A/F_A\ar[r]^{\Psi\,\,\,\,\,\,}\ar[d]& A(G/F_G)\ar[d] \\
\bar M\ar[r]_{\id_M}&M/F_M}
\end{xy}
\end{align*}
by 
\[\Psi([a_m])=A(\pi_G)(a_m)\]
for all $a_m\in A$. To see that this does not depend on the representative, note that
 $F_\core=\ker(A(\pi_G))$ and
recall that $a_m\sim a_n$ 
if and only if $\bar a_m$
is the $\nabla$-parallel transport
of $\bar a_n$
along a path lying in the leaf through $m$ of $F_M$ (Corollary \ref{holonomy_cor}). 
Without loss of generality, there exists a $\nabla$-parallel section $a\in\Gamma(A)$ such that
$a(m)=a_m$ and $a(n)=a_n+b_n$ for some $b_n\in F_\core(n)$.
Then, using \eqref{the_equation}, we get
\begin{align*}
\Psi([a_m])&=A(\pi_G)(a_m)
=\left(A(\pi_G)\circ a\right)(m)=\left(\tilde a\circ\pi_M\right)(m)\\
&=\left(\tilde a\circ\pi_M\right)(n)
=A(\pi_G)(a_n)=\Psi([a_n]).
\end{align*}
Hence,
$\Psi$ is a  well-defined vector bundle homomorphism
over the identity on $M/F_M$.
Furthermore, the considerations above show that 
for any section $[a]$ of $A/F_A$, and  corresponding $\nabla$-parallel
section $a$ of $A$, we get 
\[\Psi\circ [a]=\tilde a.
\]
 The compatibility of the Lie algebroid brackets and anchors
is then immediate by the construction of $A/F_A$, and the fact that $A(\pi_G)$
is a Lie algebroid morphism.
\end{proof}

\vspace*{0.5cm}
\section{Further remarks and applications}
\subsection{Regular Dirac structures and the kernel of the associated presymplectic groupoids}

Let $(M,\mathsf D)$ be a Dirac manifold. Recall that $(\mathsf D\to M, \pr_{TM}, [\cdot\,,\cdot])$
is then a Lie algebroid, where $ [\cdot\,,\cdot]$ is the Courant Dorfman bracket
on sections of $TM\oplus T^*M$.

Assume that the characteristic distribution
$F_M\subseteq TM$, 
defined by 
\[F_M(m)=\{ v_m\in T_mM \mid (v_m,0_m)\in\mathsf D(m)\}
\]
for all $m\in M$,
is a subbundle of $TM$. The involutivity of  
$F_M$ follows from the properties of the Dirac structure.  Set $F_\core:=F_M\oplus\{0\}\subseteq \mathsf D$.
It is easy to check that $F_\core$ is a subalgebroid of $\mathsf D$.
Define \begin{align*}
\nabla:\Gamma(F_M)\times \Gamma\left(\mathsf D/F_\core\right)&\to \Gamma\left(\mathsf D/F_\core\right)\\
\nabla_{X}\left(\bar d\,\right)&=\overline{[(X,0),d]}. 
\end{align*} 
This map is easily seen to be a well-defined, flat, partial $F_M$-connection on $\mathsf D/F_\core$,
and the verification of the fact that $(\mathsf D, F_M,F_\core, \nabla)$
is an IM-foliation on the Lie algebroid $\mathsf D\to M$
is straightforward.

We show that if  $\mathsf D\to M$ integrates to a presymplectic groupoid $(G\rr M,\omega_G)$ \cite{BuCrWeZh04,BuCaOr09},
then $(\mathsf D, F_M,F_\core, \nabla)$
integrates to the foliation $F_G=\ker\omega_G\subseteq TG$.

The map  $\rho:=\pr_{TM}:\mathsf D\to TM$ is the anchor 
of the Dirac structure $\mathsf D$ viewed as a Lie algebroid over $M$, 
and the map $\sigma:=\pr_{T^*M}:\mathsf D\to T^*M$
defines an IM-$2$-form on the Lie algebroid $\mathsf D$ (see \cite{BuCrWeZh04}, \cite{BuCaOr09}).
Note that $F_\core$ is the kernel of $\sigma$
and $F_M$ is the kernel of  $\sigma^t:TM\to\mathsf D^*$.
The  two-form  $\Lambda:=\sigma^*\omega_{\rm can}\in\Omega^2(\mathsf D)$ is 
morphic in the sense 
that 
\begin{align*}
\begin{xy}
\xymatrix{
T\mathsf D\ar[r]^{\Lambda^\sharp}\ar[d]&T^*\mathsf D\ar[d]\\
TM\ar[r]_{-\sigma^t}&\mathsf D^*
}
\end{xy}
\end{align*}
is a Lie algebroid morphism (\cite{BuCaOr09}). See, for instance \cite{Mackenzie05}, for the Lie algebroid 
structure on $T^*\mathsf D\to \mathsf D^*$. If $\mathsf D\to M$ integrates
to a presymplectic groupoid $(G\rr M, \omega_G)$, the Lie algebroid $T^*\mathsf D\to \mathsf D^*$
is isomorphic to the Lie algebroid of the cotangent groupoid $T^*G\to \mathsf D^*$
and the map $\Lambda^\sharp$ integrates via the identifications 
$T\mathsf D\simeq A(TG)$ and $T^*\mathsf D\simeq A(T^*G)$
to the vector bundle map $\omega_G^\sharp$, that is a Lie groupoid morphism. See \cite{BuCaOr09} for more details.

We show that the morphic foliation $F_{\mathsf D}\subseteq T\mathsf D$
defined by $(\mathsf D, F_M,F_\core, \nabla)$ as in Theorem \ref{main_theorem}
is equal to the kernel of $\Lambda^\sharp$.

Let $n$ be the dimension of $M$ and $k$ the rank of $F_M$. Then $\mathsf D$  
is spanned locally
by frames  of $n$ parallel sections, the first $k$ of them spanning $F_\core$.
If $d$ is a parallel section of $\mathsf D$, we have $\ldr{X}d\in\Gamma(F_\core)$
for all $X\in\Gamma(F_M)$, that is, $\ldr{X}(\sigma(d))=0$
for all $X\in\Gamma(F_M)$. Since $\ip{X}\sigma(d)=0$
for all $X\in\Gamma(F_M)$, this yields
$\ip{X}\dr(\sigma(d))=0$ for all $X\in\Gamma(F_M)$.
Hence, using this type of frames, we find with formulas (4.57) and (4.58) in 
\cite{BuCaOr09}, that the kernel of $\Lambda^\sharp$
is spanned by the restriction to $F_M$ of   the linear sections 
defined by parallel sections of $\mathsf D$, 
and by the restrictions to $F_M$ of the core sections 
defined by sections of $F_\core$. Hence, by construction, the foliation $F_{\mathsf D}$
is the kernel  of $\Lambda^\sharp$.

Since the kernel of $\omega_G^\sharp$
is multiplicative with Lie algebroid equal to the kernel of 
$\Lambda^\sharp$, this yields $F_G=\ker\omega_G^\sharp$.

Note that if $F_M\subseteq TM$ is simple, then the leaf space 
$M/F_M$ has a natural Poisson structure such that 
the projection $(M,\mathsf D)\to (M/F_M,\pi)$
is a forward Dirac map. 
Under a completeness condition and if $F_G\subseteq TG$ is also simple,
 we get a Lie groupoid  $G/F_G\rr M/F_M$, with a natural 
symplectic structure  $\omega$ such that the projection $\pi_G:G\to G/F_G$
satisfies $\pi^*_G\omega=\omega_G$. It would be interesting to study 
the relation between the integrability of the Poisson manifold 
$(M/F_M,\pi)$ and the completeness conditions on $F_G$
(see \cite{Jotz11b}) so that the quotient $(G/F_G\rr M/F_M,\omega)$
is a symplectic groupoid.

\subsection{Foliated algebroids in the sense of Vaisman}
In \cite{Vaisman10b}, foliated Lie algebroids are defined as follows. A foliated Lie algebroid
is a Lie algebroid $A\to M$ together with a subalgebroid $B$ of $A$ and an involutive 
subbundle $F_M\subseteq TM$  such that 
\begin{enumerate}
\item $\rho(B)\subseteq F_M$,
\item  $A$ is locally spanned 
over $C^\infty(M)$ by \emph{$B$-foliated cross sections}, i.e.,  sections $a$ of $A$ such that $[a,b]\in\Gamma(B)$ 
for all $b\in B$.
\end{enumerate}

Recall our definition of IM-foliation on a Lie algebroid (Definition \ref{def_foliated_algebroid}).
Since the $F_M$-partial connection is flat, we get by Proposition \ref{parallel_sections}
the existence of frames of parallel sections for $A$. By the properties of the connection, these
are  \emph{$F_\core$-foliated cross sections}. Since $(1)$ is also 
satisfied by hypothesis, our IM-foliations
are foliated Lie algebroids in the sense of Vaisman if we set $B:=F_\core$.

The object that integrates the foliated algebroid in the sense of Vaisman is
the right invariant image of $B$, which defines a foliation on $G$ that is tangent to the $\s$-fibers
and invariant under left multiplication. This is exactly 
the intersection of our multiplicative subbundle $F_G\subseteq TG$
integrating $(A\to M, F_M, F_\core,\nabla)$ with $T^\s G$.

\vspace*{0.5cm}

\appendix
\section{Invariance of bundles under flows}
We  prove here a result that is standard, but the proof of which  
is difficult to find in the literature.
\begin{theorem}\label{appendix_thm}
Let $M$ be a smooth manifold and $E$ be a subbundle of the direct sum vector bundle $\mathbb{T}M:=TM\oplus T^*M$.
Let $Z\in\mx(M)$ be a smooth vector field on $M$ and denote its flow by $\phi_t$. 
If 
\[ \ldr{Z}e\in\Gamma(E)\quad \text{ for all }
\quad e\in\Gamma(E),\]
then 
\[\phi_t^*e\in\Gamma(E)\quad \text{ for all }
\quad e\in\Gamma(E) \quad \text{ and } \quad t\in \R \,\, \text{where this makes sense.}
\]
\end{theorem}

\begin{corollary}\label{appendix_cor}
Let $F$ be a subbundle of the tangent bundle $TM$ of  a smooth manifold $M$.
Let $Z\in\mx(M)$ be a smooth vector field on $M$ and denote its flow by $\phi_t$. 
If 
\[[Z,\Gamma(F)]\subseteq \Gamma(F),\]
then 
\[T_m\phi_t F(m)=F(\phi_t(m))
\]
for all $m\in M$ and $t$ where this makes sense.
\end{corollary}

\begin{proof}
Choose $X\in\Gamma(F)$ and $m\in M$. Then, by Theorem \ref{appendix_thm}, we have 
$T\phi_{t}\circ X\circ\phi_{-t}=\phi_{-t}^*X\in\Gamma(F)$
for all $t$ where this makes sense, and hence:
$$T_m\phi_t X(m)=\left(\phi_{-t}^*X\right)(\phi_t(m))\in F(\phi_t(m)).$$
\end{proof}

\begin{proof}[Proof of Theorem \ref{appendix_thm}]
The subbundle $E$ of $\mathbb{T}M$ is an embedded submanifold of $\mathbb{T}M$.
For each section $\sigma$ of $\mathbb{T}M$, the smooth function 
$l_\sigma:\mathbb{T}M\to \R$ is defined by 
$$l_\sigma(v,\alpha)=\langle\sigma(p(v,\alpha)),(v,\alpha)\rangle$$ 
for all $(v,\alpha)\in \mathbb{T}M$, where $p:\mathbb{T}M\to M$ is the projection.
For all $e\in E$, the tangent space $T_e E$ of the submanifold $E$ of $\mathbb{T}M$ is 
equal to $$\ker\left\{\dr_el_\sigma\mid \sigma\in\Gamma\left(E^\perp\right)\right\}.$$

Consider the complete lift $\tilde Z$ to $\mathbb{T}M$ of $Z$, i.e.,
the vector field $\tilde Z\in\mx(\mathbb{T}M)$ defined 
by
\[ \tilde Z(l_\sigma)=l_{\ldr{Z}\sigma}\qquad \text{ and } \qquad
\tilde Z(p^*f)=p^*(Z(f))\]
for all $\sigma\in\Gamma(\mathbb{T}M)$ and $f\in C^\infty(M)$ (see \cite{Mackenzie05}).

Choose $e\in E$ and $\sigma\in\Gamma\left(E^\perp\right)$. Then we
have $\ldr{Z}\sigma\in\Gamma\left(E^\perp\right)$
since for all $\tau\in\Gamma(E)$:
\[\left\langle\ldr{Z}\sigma, \tau\right\rangle=Z\left(\langle\sigma,\tau\rangle\right)-
\left\langle\sigma,\ldr{Z}\tau\right\rangle=0.\]
This leads to \[(\dr_el_\sigma)(\tilde Z(e))
=\left(\tilde Z(l_\sigma)\right)(e)
=l_{\ldr{Z}\sigma}(e)=0.
\]
Hence, the vector field $\tilde Z$ is tangent to $E$ on $E$. As a consequence,
its flow curves starting at points of $e$
remain in the submanifold $E$.

It is easy to check that the flow $\Phi_t$
of the vector field $\tilde Z$ is equal to $(T\phi_t,(\phi_{-t})^*)$,
i.e.,
$$\Phi_t(v_m,\alpha_m)=(T_m\phi_t(v_m),\alpha_m\circ T_{\phi_t(m)}\phi_{-t})$$
for all $(v_m,\alpha_m)\in \mathbb{T}M(m)$.
Choose a section $(X,\alpha)\in\Gamma(E)$ and a point $m\in M$. We find 
\begin{align*}
(\phi_t^*(X,\alpha))(m)
&=\left(T_{\phi_t(m)}\phi_{-t}X(\phi_t(m)), \alpha_{\phi_t(m)}\circ T_m\phi_t\right)
=\Phi_{-t}\left((X,\alpha)(\phi_t(m))\right)
\in  E(m)
\end{align*}
since $(X,\alpha)(\phi_t(m))\in E(\phi_t(m))$. Thus, we have shown 
that $\phi_t^*(X,\alpha)$ is a section of $E$.
\end{proof}

\bigskip

Assume now that $q_A:A\to M$ is a vector bundle, 
and consider a linear vector field $X$
on $A$, i.e., the map 
$X:A\to TA$ is a vector bundle homomorphism
over $\bar X:M\to TM$ such that $X\sim_{q_A}\bar X$.
Let $\phi^X_\cdot$ be the flow of $X$ and $\phi^{\bar X}_\cdot$ the
flow of $\bar X$. Then $\phi^X_t:A\to A$ is a vector bundle homomorphism
over $\phi^{\bar X}_t$ for all $t\in\R$ where this is defined. 

Recall that for any $a\in\Gamma(A)$, 
the section $D_Xa\in\Gamma(A)$ is defined by
\[(D_Xa)(m)=\left.\frac{d}{dt}\right\an{t=0}\phi^X_{-t}(a(\phi^{\bar X}_t(m)))\]
for all $m\in M$.
In the same manner, if $\varphi\in\Gamma(A^*)$, we can define
\[(D_X\varphi)(m)=\left.\frac{d}{dt}\right\an{t=0}(\phi^X_t)^*(\varphi(\phi^{\bar X}_{t}(m)))\]
for all $m\in M$.
We have then $\varphi(a)\in C^\infty(M)$, 
and 
\begin{equation}\label{bla}
\bar X(m)(\varphi(a))
=\varphi(D_Xa)(m)+(D_X\varphi)(a)(m).
\end{equation}
We can hence show the following lemma.

\begin{lemma}\label{appendix_util}
Let $A$ be a vector bundle and  $B\subseteq A$  a subbundle. 
\begin{enumerate}
\item If $(X,\bar X)$ is a linear vector field on $A$ such that
\[D_Xb\in\Gamma(B)\]
for all $b\in\Gamma(B)$, then 
$\phi^X_t(b_m)\in B\left(\phi^{\bar X}_t(m)\right)$
for all $b_m\in B_m$. 
\item Assume furthermore that $a\in\Gamma(A)$ is such that $a(m)$ is linearly independent to $B(m)$
for all $m$ in $\dom(a)$ and 
\[D_Xa\in\Gamma(B).\]
Then 
\[\phi^X_t(a(m))\in a\left(\phi^{\bar X}_t(m)\right)+B\left(\phi^{\bar X}_t(m)\right)\]
for all $m\in U$ and $t\in\R$ where this makes sense.
\end{enumerate}
\end{lemma}

\begin{proof}
\begin{enumerate}
\item
We check that the vector field $X$ is tangent to $B$ on points in $B$. Let $\varphi\in\Gamma(A^*)$
be a section that vanishes on $B$, i.e., $\varphi_m(b_m)=0$
for all $b_m\in B$. Let $l_\varphi\in C^\infty(A)$ 
be the linear function defined by $\varphi$. By \eqref{bla}, we have 
then $D_X\varphi\in\Gamma(B^\circ)$.
Choose $b_m\in B$. We have then  
\begin{align*}
\dr_{b_m}l_\varphi(X(b_m))
&=\left.\frac{d}{dt}\right\an{t=0}l_\varphi(\phi_t^X(b_m))
=\left.\frac{d}{dt}\right\an{t=0}\varphi_{\phi_t^{\bar X}(m)}(\phi_t^X(b_m))\\
&=(D_X\varphi)(b_m)=0.
\end{align*}
Thus, $X$ is tangent to $B$ on $B$ and the flow of $X$ preserves $B$.

\item Assume now that $(b_1,\ldots,b_k)$ is a local frame for $B$ on an open set $U\subseteq M$.
Complete this frame to a local frame $(b_1,\ldots,b_n)$ 
for $A$ defined on an open $U$ such that 
 $b_{k+1}:=a\in\Gamma(A)$. Let $\varphi_1,\ldots,\varphi_{n}$ be a frame for
$A^*$ that is dual to  $(b_1,\ldots,b_n)$, i.e., such that $(\varphi_{k+1},\ldots,\varphi_n)$
is a frame for  $B^\circ$ and  
$\varphi_{k+1}(a)=1$.
Then, the closed submanifold $C$ of $A\an{U}$
defined by $C(m)=a(m)+B(m)$
is the level set with value $(1,0,\ldots,0)$ of the function 
\[(l_{\varphi_{k+1}},\ldots,l_{\varphi_n}):A\an{U}\to \R^{n-k}.
\]
Since $D_Xa\in\Gamma(B)$ for the linear vector field
$(X,\bar X)$ on $A$,
we get 
\[0=\bar X(\varphi_i(a))=\varphi_i(D_Xa)+D_X\varphi_i(a)=0+D_X\varphi_i(a)\]
for $i=k+1,\ldots,n$
and this yields
as before for all $b_m\in B$:
\begin{align*}
\dr_{a(m)+b_m}l_{\varphi_i}(X(a(m)+b_m))
&=\left.\frac{d}{dt}\right\an{t=0}l_{\varphi_i}(\phi_t^X(a(m)+b_m))\\
&=(D_X\varphi_i)(a(m)+b_m)=0.
\end{align*}
Hence, $X$ is tangent to $C$ on points of $C$. That is, the flow of $X$
preserves $C$ and the proof is finished.
\end{enumerate}
\end{proof}

\bibliographystyle{amsplain}

\def\cprime{$'$} \def\polhk#1{\setbox0=\hbox{#1}{\ooalign{\hidewidth
  \lower1.5ex\hbox{`}\hidewidth\crcr\unhbox0}}}
\providecommand{\bysame}{\leavevmode\hbox to3em{\hrulefill}\thinspace}
\providecommand{\MR}{\relax\ifhmode\unskip\space\fi MR }
\providecommand{\MRhref}[2]{%
  \href{http://www.ams.org/mathscinet-getitem?mr=#1}{#2}
}
\providecommand{\href}[2]{#2}

\end{document}